\documentclass[11pt,a4paper]{article}
\usepackage{amsmath,amsthm,amsopn,amsfonts,graphicx,amssymb,dsfont,color,yfonts}
\usepackage{mathrsfs}
\usepackage{color}
\usepackage{geometry}
\usepackage{subfigure}
\usepackage{algorithm}
\usepackage{algpseudocode}
\newcommand{\ds}{\displaystyle}
\newcommand{\mb}{\overline{ m}}
\newcommand{\rb}{\overline{ r}}
\newcommand{\xb}{\overline{ x}}

\newcommand{\ts}{v}
\newcommand{\tz}{\tilde z}
\newcommand{\tts}{\beta}

\newcommand{\bu}{\mathbf{u}}

\newcommand{\bal}{\mathbf{a}}
\newcommand{\m}{\textswab{M}}
\newcommand{\mm}{\textswab{m}}
\newcommand{\kk}{\mathbf{k}}
\newcommand{\x}{\mathbf{x}}

\newcommand{\ty}{\tilde{y}_1}
\newcommand{\Oc}{\mathcal  O}

\newcommand{\p}{\partial}
\newcommand{\hy}{\hat{y_1}}
\newcommand{\hm}{\hat{m}}

\newcommand{\abs}[1]{\left\lvert#1\right\rvert}
\newcommand{\norm}[1]{\left\lVert#1\right\rVert}

\renewcommand{\phi}{\varphi}

\def\lp {\left( }
\def\rp {\right) }
\def\R{\mathbb R}
\def\N{\mathbb N}
\def\C{\mathcal C}
\def\M{\mathcal M}

\definecolor{aquamarine}{rgb}{0.13, 0.68, 0.8} % couleur de Yoann
\definecolor{Lionel}{rgb}{0.9, 0.5, 0.15}
\def\Li{\color{black}}

\def\Bk{\color{black}}
\def\Rd{\color{black}}

\newtheorem{theorem}{Theorem}[section]
\newtheorem{proposition}[theorem]{Proposition}
\newtheorem{corollary}[theorem]{Corollary}

\theoremstyle{definition}

\newtheorem{lem}[theorem]{Lemma}

\title{Adaptation in general temporally changing environments}

\author{L. Roques$^{\hbox{\small{ a}}}$,  F. Patout$^{\hbox{\small{ a}}}$,  O. Bonnefon$^{\hbox{\small{ a}}}$,  and G. Martin$^{\hbox{\small{ b}}}$  \\
\\
\footnotesize{$^{\hbox{a }}$INRAE, BioSP, 84914, Avignon, France}\\
\footnotesize{$^{\hbox{b }}$ ISEM (UMR 5554), CNRS, 34095, Montpellier, France}
}

\date{}

\begin{document}

\maketitle

\begin{abstract}
We analyze a nonlocal PDE model describing the dynamics of adaptation of a phenotypically structured population, under the effects of mutation and selection, in a changing environment. Previous studies have analyzed the large-time behavior of such models, with particular forms of environmental changes, either linearly changing or periodically fluctuating.
We use here a completely different mathematical approach, which allows us to consider very general forms of environmental variations and to give an analytic description of the full trajectories of adaptation, including the transient phase, before a stationary behavior is reached. The main idea behind our approach is to study a bivariate distribution of two `fitness components' which contains enough information to describe the distribution of fitness at any time. This distribution solves a degenerate parabolic equation that is dealt with by defining a multidimensional cumulant generating function associated with the distribution, and solving the associated transport equation.

We apply our results to several examples, and check their accuracy, using stochastic individual-based simulations as a benchmark. These examples illustrate the importance of being able to describe the transient dynamics of adaptation to understand the development of drug resistance in pathogens.

\end{abstract}

\section{Introduction and main assumptions}

Understanding the impact of external factors on the dynamics of fitness distributions in asexuals is a fundamental issue in population genetics, with implications for the evolution of microbial pathogens such as viruses, bacteria and cancer cells. Drug resistance may occur when a pathogenic organism (\textit{e.g.}, a bacteria in presence of an antibiotic) manages to reach a positive growth rate (equivalently absolute fitness) due to genetic adaptation. Being able to describe the effect of various types of environmental changes on the trajectories of adaptation is therefore a crucial issue for the elaboration of drug resistance management strategies, to which mathematical models may help answer.

Recent models of asexual adaptation based on partial differential equations (PDEs) or  integro-differential equations (IDEs)
typically describe the dynamics of the distribution of a single phenotypic trait in a fixed environment. This trait can be fitness itself as in \cite{AlfCar14,GilHamMarRoq17,GilHamMarRoq19,TsiLev96}, or a given trait $x \in \R$ determining fitness, as in \cite{AlfCar17,AlfVer18,ChaFer13,HamLav19}, leading to equations of the form:
$$\partial_t q(t,x)=\mathcal{M}[t,x,q(t,x)]+q(t,x)\, (m(x)-\mb(t)).$$Here, $\mathcal{M}$ is a differential or an integral operator describing the effect of mutations on the distribution $q(t,\cdot)$ of the trait $x$. The last term $ q(t,x)\,(m(x)-\mb(t))$ corresponds to the effects of selection, see \textit{e.g.}~\cite{TsiLev96}: $m(x)$ is a function which describes the relationship between the trait $x$ and fitness, and $\mb(t)$ is the mean fitness in the population at time~$t$. The fitness that we consider in this work is a `relative fitness'. It is connected to the Malthusian growth rate $r$ via the formula: $r(x)=r_{max}+m(x)$ ($m(x)\le 0$ and $r_{max}>0$ is a constant corresponding to the growth rate of an optimum phenotype).

In Fisher's geometrical model (FGM), a multivariate phenotype at a set of $n$  traits (a vector $\x \in \R^n$) determines fitness. The most widely used version assumes a quadratic form of the Malthusian fitness function $m(\x)$, which decreases away from a single optimum $\Oc_0 \in \R^n$, \cite{MarLen15,Ten14}:
\begin{equation} \label{eq:m(x)}
  m(\x)=-\frac{\|\x-\Oc_0\|^2}{2},
\end{equation}
with $\|\cdot\|$ the Euclidian norm in $\R^n$. To describe the mutation effects on phenotypes, the standard  `isotropic Gaussian FGM' uses a normal distribution $\mathcal N (0, \lambda \, I_n)$ with $\lambda>0$ the phenotypic mutational variance at each trait and $I_n$ the identity matrix \cite{Kim65,Lan80}.
Overall, assuming a constant mutation rate $U$ per capita per unit time, the corresponding integro-differential equation describing the dynamics of the phenotype distribution $q(t,  \x)$, under the combined effects of selection and mutation, is $\partial_t q (t, \x) = U\, \lp J\star q - q \rp +  q (t, \x)\,(m (\x)-\mb(t)), \  t >0, \  \x \in \R^n,$
with $\mb(t)$ the mean fitness in the population at time $t$,
and $\star$ the standard convolution product in $\R^n$ and $J$ the (Gaussian) probability density function associated with the normal distribution $\mathcal N (0, \lambda \, I_n)$.

In this work, we focus on the case of a changing environment: we assume that, due to an external factor (\textit{e.g.}, a drug dose, a temperature, \textit{etc}), the phenotype to fitness relationship \eqref{eq:m(x)} is changed. We take this change into account through a moving optimum, \textit{i.e.} we assume that
\begin{equation} \label{eq:m(t,x)}
  m(t,\x)=-\frac{\|\x-\Oc(t)\|^2}{\Li 2 \Bk},
\end{equation}
with
\begin{equation} \label{eq:O(t)}
  \Oc(t)=\Oc_0+ \delta(t) \, \bu ,
\end{equation}
with $\delta(t) \in \C(\R_+)$, $\delta(0)=0$ and  $\bu$ a unit vector in $\R^n$ (without loss of generality, we assume in the sequel that $\Oc_0=0$, and $\bu=(1,0,\ldots,0)$). In such case, the equation describing the dynamics of the phenotype distribution becomes:
\begin{equation}
\partial_t q (t, \x) = U\, \lp J\star q - q \rp + q(t, \x)\,  (m (t,\x)-\mb(t)), \  t >0, \  \x \in \R^n,
\label{eq:main_IDE_moving_env}
\end{equation}
with this time:
\begin{equation}
\label{def:mbar_t}
    \mb(t) = \int_{\R^n} m(t,\x) \, q(t, \x) \, d\x.
\end{equation}

We approach the mutational effects $U\, (J\star q -q)$ by a diffusion (Laplace) operator, leading to the main equation that is studied in this paper:
\begin{equation}
\partial_t q (t, \x) = \frac{\mu^2}{2} \Delta q+  q(t, \x) \, (m (t,\x)-\mb(t)), \  t >0, \  \x \in \R^n,
\label{eq:main_PDE_moving_env}
\end{equation}
with $\mu =\sqrt{U\, \lambda}>0$ the mutation parameter; we refer to \cite{HamLav19} (Appendix) for further details on the derivation of this diffusion approximation. The regime where it applies corresponds to the `Weak Selection Strong Mutation' (WSSM) regime, where a wide diversity of lineages accumulate mutations and co-segregate at all times.

%In our case, we work on the phenotype distribution (by definition, of integral $1$)

The main goal of our work is to describe the dynamics of the mean fitness $\mb(t)$ in the population for very general scenarios of environmental changes, i.e., with a general form for $\delta(t)$. The value of $\mb(t)$ is fundamentally connected with the question of drug resistance, or evolutionary rescue in a broader context \cite{GomHol95}, as they occur when the mean growth rate (or equivalently mean absolute fitness) in the population $\overline{r}(t):=r_{max}+\mb(t)$ becomes positive.

Several particular forms of environmental changes have already been considered. The case of an optimum shifting with a constant speed $\delta(t)=c\, t$ has inspired several developments. First, in models without adaptation, where a favorable region moves at a constant speed, as in the Fisher-KPP reaction-diffusion equations studied in~\cite{BerDie09,BerFan18,BerRos08}. Then, \cite{AlfBer17} considered again an optimum shifting with a constant speed, in a model including both a 1D space variable and a 1D phenotypic trait. From a mathematical viewpoint, this corresponds to an equation of the form~\eqref{eq:main_PDE_moving_env} in $\R^2$. The case of periodically fluctuating environments has also attracted much interest from mathematicians. In the 1D case, \cite{LorChi15} derived explicit Gaussian solutions of a PDE model describing the dynamics of a phenotype distribution with a periodically varying phenotype optimum. In \cite{CarNad20}, comparable models were considered, with phenotypes in some bounded subset of $\R^n$; the authors give conditions for the persistence of the population, based on the sign of the principal eigenvalue of a time-periodic parabolic operator, and study the large-time behavior of the solution. In \cite{FigMir18}, fitness functions $m(t,\x)$, periodic with respect to $t$ and with $\x\in \R^n$ have also been considered (see also \cite{FigMir19}). The works \cite{FigMir18,FigMir19} are based on the method of constrained Hamilton-Jacobi equations, which has been developed to study the evolution of phenotypically structured  populations, with integral or differential mutation operators (\textit{e.g.}, \cite{BarMir09,DieJab05,GanMir17,LorMir11,PerBar08}). This method assumes a small mutation parameter of order $\varepsilon \ll 1$, and is based on a scaling $t\to t/\varepsilon$. Thus, it typically describes asymptotic evolutionary dynamics, at large times and in a  `small mutation' regime. To the best of our knowledge, it cannot lead to explicit transient trajectories of adaptation. Temporally piecewise constant environments have also been studied with the same type of methods in the recent work \cite{CosEtc19}. Note that the equations that were studied in \cite{AlfBer17,CarNad20,FigMir18,FigMir19,LorChi15} have the general form:
\begin{equation}\label{eq:main_form_n}
\partial_t n (t, \x) = \frac{\mu^2}{2} \Delta n+  n(t, \x) \, (r (t,\x)-\rho(t)), \  t >0, \  \x \in \Omega \subseteq \R^n,
\end{equation}
with $n(t,\x)$ the total population density and $\rho(t)$ its integral over $\Omega$. The study of this equation is in fact equivalent to our problem~\eqref{eq:main_PDE_moving_env}: it is easily checked that $q(t,\x)=n(t,\x)/\rho(t)$ satisfies \eqref{eq:main_PDE_moving_env} with $r (t,\x)=r_{max}+m(t,\x)$, see Appendix~A.

Compared to the above-mentioned works,  we use here a completely different approach, which allows us: (i) to consider very general forms of environmental variations; (ii) to give an analytic description of the full trajectories of adaptation, including the transient phase, before a stationary behavior is reached. Our results are valid in any dimension $n$, and do not assume that the solution has a Gaussian form. The main ideas behind our approach is to study a bivariate distribution $p(t,m_1,m_2)$ of two `fitness components' which contains enough information to describe the distribution of fitness at any time $t$. The distribution $p$ solves a degenerate parabolic equation that is dealt with by defining a multidimensional cumulant generating function associated with the distribution, and solving the associated transport equation.

Our main results are presented in the next section. We begin in Section~\ref{sec:exist} with a preliminary standard existence and uniqueness results of the solution $q(t,\x)$ of the Cauchy problem associated with~\eqref{eq:main_PDE_moving_env}; then, in Section~\ref{sec:cumul} we study the distribution $p(t,\cdot)$ and derive the equation solved by the cumulant generating function; in Section~\ref{sec:galeformula}, we present our main results on the dynamics of the mean fitness in a general setting; in Section~\ref{sec:examples}, we apply these results to particular forms of the function $\delta(t)$, and we compare our results with the existing literature. In Section~\ref{sec:numm}, we compare our theoretical results with numerical simulations of a stochastic individual-based model. These sections are followed by a discussion. Proofs are presented in Section~\ref{sec:proofs}.

\section{Main results \label{sec:results}}

\subsection{Existence and uniqueness of the solution of the Cauchy problem \label{sec:exist}} The existence and uniqueness of the solution $q(t,\x)$ of \eqref{eq:main_PDE_moving_env} with initial condition $q_0$ does not follow from standard parabolic theory as the function $m(t,\x)$ is unbounded. However, they can easily be adapted from the results in \cite{HamLav19}, in order to take into account the time-dependence of $m(t,\x)$. We recall here the main arguments that lead to these existence and uniqueness results.

We need the following assumptions on the initial distribution $q_0$:
\begin{equation} \label{eq:hypq0a}
    q_0\in C^{2+\alpha}(\R^n),
\end{equation}
for some $\alpha\in (0,1)$, that is, $\|q_0\|_{C^{2+\alpha}(\R^n)}<+\infty$. Moreover, as $q_0$ is a (probability) distribution, we assume that:
\begin{equation} \label{eq:hypq0_int1}
    q_0\ge 0 \hbox{ and } \int_{\R^n}q_0(\x) d\x=1.
\end{equation}
We also assume that $q_0$ has a fast decay rate as $\|\x\|\to +\infty$, in the sense that there exists a non-increasing function $g\in \C(\R_+,\R_+)$ (with $\R_+=[0,+\infty)$) such that:
\begin{align} \label{eq:hypq0b}
0\le q_0\le g(\|\cdot\|)\hbox{ in }\R^n,
\hbox{ and } \int_{\R^n}e^{b \|\x\|}\,g(\|\x\|)\,d\x<+\infty \hbox{ for all }b>0.
\end{align}
We first recall a standard existence and uniqueness result for linear parabolic equations with unbounded coefficients in $\R^n$.
\begin{theorem}[\cite{AroBes67,Cha70}]\label{th:existence_v}
The problem
\begin{equation}
\left\{ \begin{array}{rcl}
     \partial_t v(t,  \x) & = & \ds \frac{\mu^2}{2} \Delta v +m(t,\x)\,v(t,  \x), \  t\ge0, \ \x \in  \R^n,  \\
     v(0, \x ) & = & q_0(\x), \   \x \in \R^n,
\end{array}\right.
\label{eq:v_sans_mbar}
\end{equation}
admits a unique positive bounded solution $v\in C^{1,2}(\R_+ \times \R^n)$.
\end{theorem}
Moreover, it follows from the same arguments as those in lemma~4.2 of \cite{HamLav19} that, with $v$ defined in Theorem~\ref{th:existence_v},
$$
t\mapsto \mb_v(t):=\int_{\R^n}m(t,\x)\,v(t,\x)\,d\x,
$$
is real-valued and continuous in $\R_+$ and, for every $t\ge0$, there holds:
$$
1+\int_0^t\overline{m}_v(s)\,ds=1+\int_0^t \int_{\R^n} m(t,\x) v(s,\x)\,d\x \, ds=\int_{\R^n}v(t,\x)\,d\x>0.
$$
This allows us to define:
\begin{equation}
q(t,\x) = {v(t,\x) \over 1+\int_0^t\overline{m}_v(s)\,ds},
\label{eq:rel_v_q}
\end{equation}
for every $(t,\x)\in\R_+\times\R^n$. Arguing as in theorem~4.1 of \cite{HamLav19}, it is straightforward to check that $q$ is the unique solution of \eqref{eq:main_PDE_moving_env}. More precisely,
\begin{theorem}[\cite{HamLav19}]
There exists a unique nonnegative solution $q \in C^{1,2}(\R_+ \times \R^n) $ of \eqref{eq:main_PDE_moving_env}
such that $q\in L^\infty((0,T)\times\R^n)$ for all $T>0$, and the function:
$$t\mapsto\mb(t)=\int_{\R^n}m(t,\x)\,q(t,\x)\,d\x,$$
is real-valued and continuous in $\R_+$. Moreover, we have:
$$\ds\forall\,t\geq 0,\ \int_{\R^n} q(t, \x)\, d\x = 1.$$
\label{thm:existuniq}
\end{theorem}
Additionally, as $m(t,\x)\le 0$ and $\mb(t)$ is bounded, a standard comparison argument implies that:
$$0\le q(t,\x)\le B(t) \, K\star q_0(\x),$$for some positive bounded function $B(t)$, and with $K$ the heat kernel in dimension $n$:
$$K(\x)=\frac{1}{(2\, \pi \, t \, \mu^2)^{n/2}}e^{-\frac{\|\x^2\|}{2 \,t \, \mu^2}}.$$With the assumption \eqref{eq:hypq0b}, this implies that $q$ is exponentially bounded at all times:
\begin{equation}\label{eq:q_exp_bound}
\hbox{for all }t\ge 0     \hbox{ for all }b>0, \ \int_{\R^n}e^{b \|\x\|}\, q(t,\x)\,d\x<+\infty.
\end{equation}

\subsection{Fitness components and cumulant generating functions \label{sec:cumul}}
In the previous work \cite{HamLav19}, where the optimum $\Oc_0$ remained constant, it was shown that the distribution of fitness, say $p(t,m)$, satisfies a 1D degenerate parabolic PDE. Defining the cumulant generating function$$C(t,z)=\ln \lp \int_{\R}p(t,s)\, e^{s \, z} \, ds \rp,$$ associated with this distribution, an analytically tractable 1D transport equation for $C(t,z)$ was obtained, leading to an explicit formula for $\mb(t)=\partial_z C(t,0)$ (which was consistent with the formula in \cite{MarRoq16} in the isotropic case).

Here, due to the time-dependence of $\Oc(t)$, one cannot expect to obtain a single autonomous PDE for the fitness distribution $p(t,m)$. Consider for instance two initial distributions of $q(0,\x)$ which are symmetric with respect to $\Oc_0=0$: $q_1(0,\x)=q_2(0,-\x)$; then the corresponding initial distribution of fitness is the same, as $m(0,\x)=m(0,-\x)$:  $p_1(0,m)=p_2(0,m)$. However, it is natural to expect that $p_1(t,m)\not \equiv p_2(t,m)$ for $t>0$: if $q_1(0,\cdot)$ is localised around the position of the optimum at later times, $q_2(0,\cdot)$ is localised in the opposite direction.

\begin{figure}
\center
\includegraphics[width=0.59\textwidth]{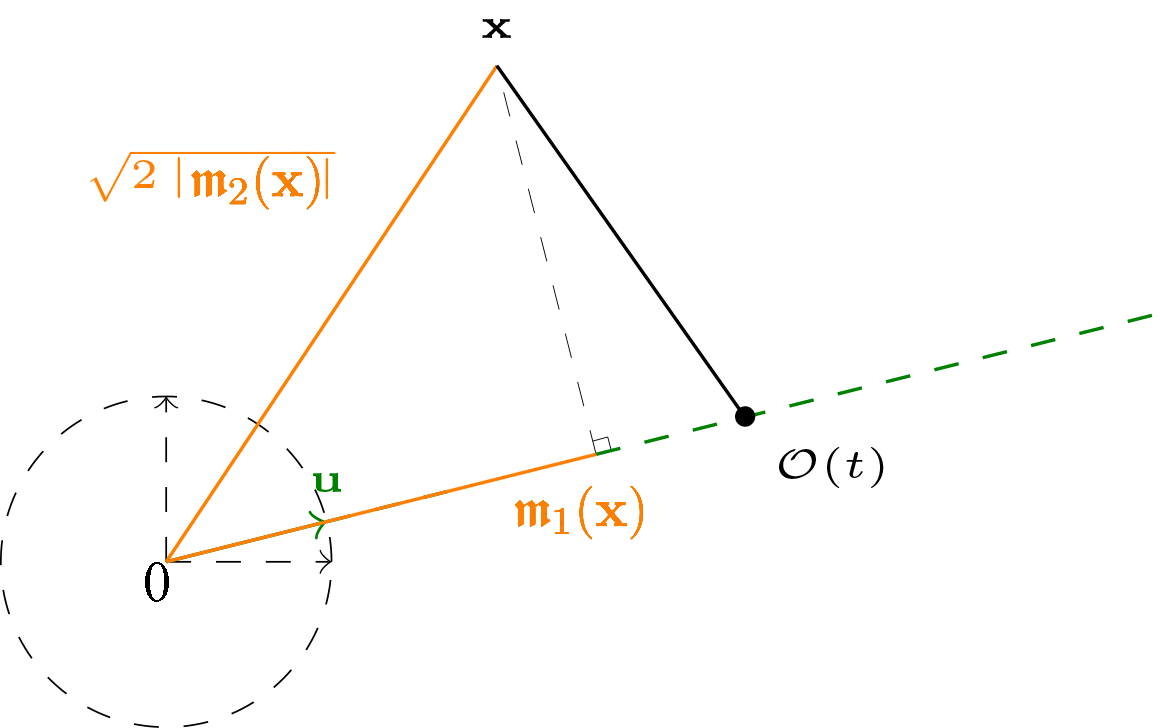}
\caption{{\bf Schematic illustration of the fitness components.} Given $\Oc(t)$, the fitness of a phenotype $\x$, $m(t,\x)=-\|\x-\Oc(t)\|^2/2$, is uniquely defined by the couple $(\mm_1(\x),\mm_2(\x))$.
}\label{fig:comp fitness}
\end{figure}

Thus, instead of focusing on the fitness distribution, we define two time-independent `components':
\begin{equation}\label{components fitness}
    \left\{
    \begin{array}{rl}
         \mm_1(\x) &= \bu \cdot \x=x_1,  \\
         \ds \mm_2 (\x)& \ds = -\frac{\|\x\|^2}{2},
    \end{array}
    \right.
\end{equation}
and we denote $\m(\x)=(\mm_1(\x),\mm_2(\x))$ for all $\x \in \R^n$, see Fig.~\ref{fig:comp fitness}. We observe that, at any time $t$, the fitness associated with $\x$ is uniquely determined from its components by the following formula:
\begin{equation} \label{eq:m(t,x)_b}
  m(t,\x)=   \ds  -\frac{\|\x-\Oc(t)\|^2}{2 }=\ds \delta(t)\, \mm_1(\x)+\mm_2(\x)-\frac{\delta(t)^2}{2}.
\end{equation}

We define $p(t,m_1,m_2)$ the (bivariate) distribution of the components $(m_1,m_2)$ at time $t$. More precisely, $p$ is defined in the next theorem.
\begin{theorem}\label{theo equation p}
There exists a unique nonnegative density function $p\in C^1(\R_+,L^2(\R \times \R_-))$ that satisfies the following relationship \begin{align}\label{link q and p}
\int_{\R^n} q(t,\x) \phi(\m(\x)) d\x = \int_{\R \times \R_-} p(t,m_1,m_2) \phi(m_1,m_2)dm_1dm_2,
\end{align} for every test functions $\phi \in L^2(\R \times \R_-)$ and all $t\ge 0$.
\end{theorem}
Using \eqref{eq:m(t,x)_b}, we observe that the mean fitness in the population at time $t$ is given by

\begin{equation}\label{eq:mbar_defp}
\mb(t)=     \ds \int_{\R^n}m(t,\x)\,q(t,\x)\,d\x
     =  \ds\int_{\R \times \R_-} p(t,m_1,m_2)  \lp \delta(t)\, m_1 + m_2 -\frac{\delta(t)^2}{2} \rp dm_1 dm_2.
\end{equation}
Similarly, the fitness variance in the population,
\begin{equation} \label{def:variance}
 V_m(t):= \int_{\R^n} m(t,\x)^2 \,  q(t,\x) \, d\x -\mb(t)^2,
\end{equation}
is given by
\begin{align}\label{eq:Vm_defp}
    V_m(t)= \delta(t)^2 V_{m_1}(t) +  V_{m_2}(t) + 2 \delta(t) cov_m(t),
\end{align}
with, for $j=1,2$    $\ds V_{m_j}(t)  = \int_{\R \times \R_{-}} p(t,m_1,m_2) m_j^2  dm_1 dm_2- \lp \int_{\R \times \R_{-}} p(t,m_1,m_2) m_j dm_1 dm_2  \rp^2,$
and:
\begin{multline*}
  cov_m(t) =
 \int_{\R \times \R_{-}} p(t,m_1,m_2) m_1 m_2  dm_1 dm_2 \\- \lp  \int_{\R \times \R_{-}} p(t,m_1,m_2) m_1  dm_1 dm_2 \rp \lp \int_{\R \times \R_{-}} p(t,m_1,m_2) m_1  dm_1dm_2 \rp.
\end{multline*}
Moreover, we can define the `cumulant generating function' associated with $p$.
\begin{theorem}\label{thm:CGF}
The cumulant generating function of the components $m_1,$ $m_2$:
\begin{equation}\label{def:C(t,z1,z2)}
C(t,z_1,z_2):=\ln \lp \int_{\R \times \R_{-}}p(t,m_1,m_2)\, e^{m_1 \, z_1 +m_2 \, z_2} \, dm_1 \, dm_2 \rp,
\end{equation}
for all $t\ge 0$, $z_1\in\R$ and $z_2\in \R_+$
 is well-defined and belongs to  $\mathcal C^{1} (\R_+\times \R \times \R_+)$. It satisfies the following equation, for $t\ge 0$ and  $(z_1,z_2)\in \R \times \R_+$:
\begin{equation}
\label{eq:C}
\left\{
\begin{array}{l}
\partial_t C(t,z_1,z_2)= \bal(t) \cdot (\nabla C(t,z_1,z_2) -\nabla C(t,0,0)) +  \kk(z_1,z_2) \cdot \nabla C(t,z_1,z_2) + \gamma(z_1,z_2), \\
C(0,z_1,z_2)=C_0(z_1,z_2), \\
C(t,0,0)=0,
\end{array}\right.
\end{equation}
\begin{equation} \label{def:gamma}
\hbox{where }\bal(t)=(\delta(t), 1)\in \R^2\hbox{ and }\left\{
\begin{array}{l}
\kk(z_1,z_2)=-\mu^2 (z_1\, z_2,z_2^2),    \\
\gamma(z_1,z_2)= \mu^2\, (z_1^2/2-n \, z_2/2).
\end{array}\right.
\end{equation}
\end{theorem}
The cumulant generating function leads to simple characterizations of the central moments of the distribution. For instance, the mean fitness and the variance can easily be computed from $C$. For $j=1,2$  $\partial_j C(t,0,0)= \int_{\R \times \R_{-}}m_j \, p(t,m_1,m_2) \, dm_1 \, dm_2,$  $\partial_{jj} C(t,0,0)= V_{m_j}(t),$ and $\partial_{1,2} C(t,0,0)= cov_m(t)$. Using \eqref{eq:mbar_defp}, and \eqref{eq:Vm_defp},
we get:
\begin{align}\label{eq:mbar_CGF}
\mb(t)& =\delta(t)\,\partial_1 C(t,0,0)+\partial_2 C(t,0,0)-\frac{\delta(t)^2}{2},\\
V_m(t) & =  \delta(t)^2   \partial_{11} C(t,0,0) +  \partial_{22} C(t,0,0) +2 \delta(t) \partial_{12} C (t,0,0). \label{eq:var_cgf}
\end{align}

\subsection{General formulas for the mean fitness and the fitness variance \label{sec:galeformula}}
\Li In order to solve the equation \eqref{eq:C} satisfied by $C$, we first note that a simpler problem can be solved explicitly. Namely, we have the following proposition.
\begin{proposition}\label{prop:solQ}
Let $\beta\in \C^1(\R_+^3,\R)$ with $\beta(t,0,0)=0$ for all $t\ge 0$, and $Q_0\in \C(\R_+^2,\R)$. For $t\ge0$ and $(z,\tz)\in\R_+^2$, the problem
\begin{equation}
\label{eq:Q_inia}
\left\{\begin{array}{l}
\partial_t Q(t,z,\tz)= (\p_z Q +\p_{\tz} Q)(t,z,\tz) -(\p_z Q +\p_{\tz} Q)(t,0,0)  + \beta(t,z,\tz), \\
Q(0,z,\tz)=Q_0(z,\tz), \\
Q(t,0,0)=0, \end{array} \right.
\end{equation}
 admits a unique solution, which is given by the expression:
\begin{equation}\label{eq:solQ} Q(t,z,\tz)= \int_0^t \beta(t-s,z+s,\tz+s)- \beta(t-s,s,s) \, ds + Q_0(z+t,\tz+t)-Q_0(t,t). \end{equation}
\end{proposition}
Next, we look for a change of variables such that the function $C$ in the rescaled variables solves a system of the form~\eqref{eq:Q_inia}. This leads to our main theorem, which can be stated as follows.
\begin{theorem}
\label{th:formule_C}
For each $t\ge 0,$ define $\phi_t: \ \R_+^2  \to  \R\times \R_+,$  by $\phi_t(z,\tz)=(y_1(t,z,\tz),y_2(z))$, with $$\left\{
\begin{array}{rl}
y_1(t,z,\tz):=& \ds\int_0^{z} \delta(z+t-s) \, \frac{\cosh(\mu \, s)}{\cosh(\mu \, z)}\, ds+(z-\tz)\, \frac{\cosh(\mu(z+t))}{\cosh(\mu \, z)},    \vspace{2mm}   \\
y_2(z):=& \ds\frac{\tanh(\mu\, z)}{\mu}.
\end{array}\right.
$$Let $\beta(t,z,\tz):=\gamma(\phi_t(z,\tz))$ with $\gamma$ defined in \eqref{def:gamma} and $Q$ be defined by~\eqref{eq:solQ}. Then, for all $t\ge 0$ and $(z,\tz)\in \R_+^2$, the solution of \eqref{eq:C} satisfies:
\begin{equation}\label{eq:C_Q}
    C(t,\phi_t(z,\tz))=Q(t,z,\tz),
\end{equation}
\end{theorem}
This result leads to an explicit expression for $\mb(t)$, as stated below.
\begin{corollary}
\label{cor:formule_mbar}
   Let $Q$ be defined as in Theorem~\ref{th:formule_C} and $R(t,z):=Q(t,z,z).$ \Bk
The mean fitness in \eqref{eq:main_PDE_moving_env} is given by
\begin{equation}\label{eq:formule_gale_mbar}
\overline{m}(t)= \partial_z R(t,0)-\delta(t)^2/2,
\end{equation}
or, more explicitly,
  \begin{equation}\label{eq:mbar_cor}
   \mb(t)=-\mu \, \frac{n}{2} \, \tanh(\mu \,t) -\frac{1}{2}\lp H_\delta(t)-\delta(t) \rp^2+  R_0'(t),
  \end{equation}
with $\ds H_\delta(t):=\mu \, \int_0^t\delta(u) \, \frac{\sinh(\mu \, u)}{\cosh(\mu \, t)}\,  du$ and \begin{equation}\label{eq:Q0primea}
    R_0'(t)=  \frac{1}{\cosh(\mu \, t)}\lp \delta(t)-H_\delta(t) \rp  \partial_1 C_0(\phi_0(t,t)) +(1- \tanh^2(\mu \,t ))\,\p_2 C_0(\phi_0(t,t)).
\end{equation}
\end{corollary} \Bk
The term $-\mu \, (n/2) \, \tanh(\mu \,t)$ in \eqref{eq:mbar_cor} corresponds to the dynamics of $\mb(t)$ with a steady optimum ($\delta\equiv 0$), and is consistent with the results in \cite{GilHamMarRoq19,HamLav19,MarRoq16}. The second term $-\frac{1}{2}\lp H_\delta(t)-\delta(t) \rp^2$ is a sort of squared distance between the position of the optimum at time $t$, and a `weighted history' of $\delta$ for $u\in (0,t)$.

The dependence of the dynamics of $\mb(t)$ with respect to the initial phenotype distribution clearly shows up in  Corollary~\ref{cor:formule_mbar}, through $R_0'(t)$. Let $C_{clonal} $ be the solution of \eqref{eq:C}, with an initial condition $C_{0,clonal}(z_1,z_2)=0$ corresponding to a clonal population at the optimum at $t=0$, \textit{i.e.} with $q_0$ a Dirac mass at $\x=0$ (the cumulant generating function associated with a Dirac mass at $\x^*$ is $z_1 \, x_1^* -z_2 \, \|\x^*\|^2/2$). Even though this initial data does not satisfy the previous conditions \eqref{eq:hypq0a}-\eqref{eq:hypq0_int1} on $q_0$, the function $C_{clonal}$ can still be defined, at least locally, through the equality \eqref{eq:C_Q}. We denote by $\mb_{clonal}(t)$ the corresponding value of the mean fitness:
$$\mb_{clonal}(t)=-\mu \, \frac{n}{2} \, \tanh(\mu \,t) -\frac{1}{2}\lp H_\delta(t)-\delta(t) \rp^2.$$Then, the mean fitness in \eqref{eq:main_PDE_moving_env} is given by
\begin{equation} \label{eq:dep_mbar_Q0}\mb(t)=\mb_{clonal}(t) +R_0'(t),\end{equation}with $R_0'$ given by~\eqref{eq:Q0primea}.

Another corollary of Theorem~\ref{th:formule_C} gives a characterization of the variance of the fitness distribution at any time.
\begin{corollary}
\label{cor:formule_variance}
The fitness variance $V_m(t)$ defined by \eqref{def:variance}  is given by
\begin{align} \label{eq:cor_variance}
V_m(t) = \p_{zz}R(t,0) + \frac{\delta'(t)}{\cosh(\mu \,t)}\partial_{\tz} Q (t,0,0). \end{align}
\end{corollary}
Other moments of the distribution of fitness could be computed as well, based on the result of Theorem~\ref{th:formule_C}. For instance, the third standardized moment of the distribution (skewness) is given by the formula \eqref{form skew} in  Appendix~B.

%The variance of fitness in \eqref{eq:mbar_CGF} is given by

\subsection{Explicit expressions for $\mb(t)$ and $V_m(\infty)$: some examples \label{sec:examples}}

%In all cases, we describe the dynamics of the mean fitness in the case of an initially clonal population, see Corollary~\ref{cor:effet_DI} and the comments above.

In this section, we apply the general formula for $\mb(t)$ derived in Corollary~\ref{cor:formule_mbar} to several particular forms of functions $\delta(t)$. Explicit but lengthy expressions for $V_m(t)$ can also be derived from Corollary~\ref{cor:formule_variance}. We only give here the formula for the asymptotic variance as $t\to +\infty$. These results are left without proof, as they are straightforward consequences of formulas~\eqref{eq:mbar_cor} and~\eqref{eq:cor_variance}.

\paragraph{Optimum shifting with a constant speed.} We make here the standard assumption (\textit{e.g.}, \cite{AlfBer17,FigMir19}) of an optimum $\Oc(t)$ which moves at a constant speed.
\begin{proposition}\label{prop:linear}
Assume that $\delta(t)=c\, t$ for some $c\in \R$. Then the mean fitness is given by:
\begin{equation}\label{eq:mbar_linear}
      \mb(t)=-\mu \, \frac{n}{2} \, \tanh(\mu \,t) -\frac{c^2}{2\, \mu^2} \tanh^2(\mu \, t)+  R_0'(t),
\end{equation}
with $$R_0'(t)= \frac{c}{\mu}\frac{\tanh(\mu \, t)}{\cosh(\mu \, t)}\p_1 C_0 (\phi_0(t,t))+(1-\tanh^2(\mu \,t ))\,\p_2 C_0(\phi_0(t,t))$$and$$\phi_0(t,t)=\lp \frac{c}{\mu^2}\lp 1-\frac{1}{\cosh(\mu \,t )}\rp, \frac{\tanh(\mu\, t)}{\mu}\rp.$$
\end{proposition}
Passing to the limit $t\to +\infty$ in \eqref{eq:mbar_linear}, we observe that $$\mb(\infty)=-\mu \frac{n}{2}-\frac{c^2}{2 \, \mu^2},$$is independent of the initial phenotype distribution. In the case $c=0$ (steady optimum), $\mb(\infty)=-\mu \, n/2$. This quantity is the `mutation load', i.e., the decrease in fitness (compared to the optimum $0$), due to mutations. When $c\neq 0$, the additional negative term $-c^2/(2 \, \mu^2)$ describes the `lag load', which corresponds here to the decrease in fitness due to the shifting in the optimum. We observe that the mutation parameter has opposite effects on the mutation and lag loads: it tends to increase the mutation load and to decrease the lag load. This leads to an optimum value $\mu^*=(2 \, c^2/n)^{1/3}$ which maximizes $\mb(\infty)$.

We recall that the growth rate of a population whose phenotype distribution satisfies~\eqref{eq:main_PDE_moving_env} can be described by $\rb(t):=r_{max}+\mb(t)$, with $r_{max}>0$ a fixed constant corresponding to the growth rate of an optimum phenotype. Persistence of the population at large times is then equivalent to $\rb(\infty):=r_{max}+\mb(\infty)>0$ (see Appendix~A). The results of Proposition~\ref{prop:linear} show that the critical speed $c^*$ for persistence ($\rb(\infty)> 0$ if $c<c^*$) or extinction ($\rb(\infty)\le 0$ if $c\ge c^*$) of the population is given by:
\begin{equation} \label{eq:maxspeed}
c^*=\mu \, \sqrt{2 \, r_{max} - \mu\, n}.
\end{equation}
Note that the condition $2 \, r_{max} - \mu\, n\ge 0$ is necessary for the survival of the population in a fixed environment.

A consequence of the results of \cite{AlfBer17} is that the critical speed satisfies $c^*=2 \sqrt{-  \lambda_\infty \, \mu^2/2}$, with $\lambda_{\infty}$
the principal eigenvalue of the operator $\phi \mapsto -\mu^2/2 \p_{xx} \phi-(r_{max} -x^2/2) \phi$ in $\R$, which is given here by $\lambda_\infty=\mu/2-r_{max}$ (principal eigenfunction: $\phi(x)=\exp(-x^2/(2\, \mu))$). Our formula \eqref{eq:maxspeed} is therefore fully consistent with the formula $c^*=2 \sqrt{-  \lambda_\infty \, \mu^2/2}$ in~\cite{AlfBer17}.

Using Corollary~\ref{cor:formule_variance}, we also obtain an explicit expression for the limit $V_m(\infty)$ of $V_m(t)$ as $t\to +\infty$:$$V_m(\infty)=\mu^2 \frac{n}{2}+\frac{c^2}{\mu}.$$Thus, the variance of the fitness distribution increases with the speed $c$: a higher speed leads to a flatter distribution. However, it is a nonmonotonic function of $\mu$: contrarily to the case of a fixed environment, $ V_m(\infty)$ first decreases with $\mu$, until a critical value which is reached at $\mu=(c^2/n)^{1/3}$, and then increases with $\mu$.

Finally using the expression for the skewness~\eqref{form skew} derived in Appendix~B, we obtain that as $t\to +\infty$, the skewness converges to:
$$\hbox{Skew}_m(\infty)=-\frac{\mu^3 \, n +3 \, c^2}{V_m(\Li\infty\Bk)^{3/2}}.$$This negative skewness implies that the distribution of fitness is asymmetrical, with a longer left tail. The skewness becomes even more negative when the speed $c$ is increased, which therefore reinforces the asymmetry of the distribution.

\paragraph{Sub- and superlinear cases.}  We assume here that the position of the optimum $\Oc(t)$ is a sublinear, or superlinear function of $t$: $\delta(t)=c\, t^\alpha,$ with $\alpha >0$, $\alpha\neq 1.$ The general formulas~\eqref{eq:formule_gale_mbar} and~\eqref{eq:cor_variance} can be applied to derive explicit expressions for $\mb(t)$ and $V_m(t)$. As these expressions are rather complex, we only summarize some asymptotic properties below.
\begin{proposition}\label{prop:nonlinear}
Assume that $\delta(t)=c\, t^\alpha$ for some $c\in \R^*$ and $\alpha>0$.

(i) If $\alpha<1$, then $\mb(t) \to -\mu \, n/2$ and $V_m(t) \to \mu^2 \, n/2$, as $t\to +\infty$.

(ii) If $\alpha>1$,  then  $\mb(t) \to -\infty$ and $V_m(t) \to +\infty$, as $t\to +\infty$.
\end{proposition}
Thus, if $\alpha<1$ (sublinear case), the lag load is equal to $0$:
at large times, the population tends to be as well adapted as in the case of a steady optimum. On the other hand, if $\alpha>1$ (superlinear case) the lag load is infinite; this means that adaptation is not possible.

\paragraph{Periodically varying optimum.} The case of an optimum $\Oc(t)$ with a periodic trajectory is particularly relevant in applications, when a population faces an external factor which is itself periodic (concentration of an antibiotic, temperature, \ldots ). We consider here a particular case, to illustrate the global shape of $\mb(t)$ in such situations: $\Oc(t)$ oscillates between the two points $\pm \delta_{max}\, \bu\in \R^n$.
\begin{proposition}\label{prop:sinus}
Assume that $\delta(t)=\delta_{max} \, \sin(\omega\, t)$
for some $\omega \in \R_+^*$ and $\delta_{max}>0$.
The mean fitness is given by
\begin{equation}\label{eq:mbar_sinus1}
      \mb(t)=-\mu \, \frac{n}{2} \, \tanh(\mu \,t) -\frac{1}{2}\lp\frac{\delta_{max} \,  \omega}{\omega^2 + \mu^2}\rp^2 ( \omega \, \sin( \omega \, t)+\mu \, \cos(
       \omega \, t) \tanh(\mu \, t))^2 +R_0'(t),
\end{equation}
and the average value of $\mb(t)$ over one period converges to
\begin{equation} \label{eq:mean_lag_sinusa}
    \langle \mb_\infty\rangle:=\lim_{t\to +\infty}  \frac{\omega}{\pi}  \int_t^{t+\pi/\omega}\mb(s) \, ds=-\mu\, \frac{n}{2}-\frac{\delta_{max}^2 \, \omega^2  }{4 \omega^2 + 4\mu^2}.
\end{equation}
\end{proposition}

Thus, asymptotically in time, the mean fitness $\mb(t)$ becomes periodic with period $\pi/\omega$. Additionally, the formula~\eqref{eq:mean_lag_sinusa} tells us that the average mean fitness at large times is a decreasing function of the frequency $\omega/\pi$: higher frequencies tend to impede adaptation. In a rapidly oscillating environment, i.e. as $\omega\to +\infty$, the average lag load converges to $-\delta_{max}^2/4$, which means that the system behaves in average as if the phenotype distribution was at a distance $\delta_{max}$ from the optimum. Conversely, in a slowly oscillating environment, i.e. as $\omega\to 0$, the average lag load is equivalent to $-\delta_{max}^2 \, \omega^2  /(4\mu^2)$; in this case, the system behaves in average as in the case of a steadily moving optimum, with speed $\delta_{max}\, \omega/\sqrt{2}$ (see Proposition~\ref{prop:linear}). Formula \eqref{eq:mean_lag_sinusa} also enables us to study the dependence of $\langle \mb_\infty\rangle$ with respect to the mutation parameter $\mu$: it is convex until the inflexion point $\mu=\omega/\sqrt{3}$ and then concave. At this inflexion point, $\p_{\mu}\langle \mb_\infty\rangle(\mu=\omega/\sqrt{3})=-n/2+3 \delta_{max}^2 \sqrt{3}/(32\, \omega)$. Thus, if $-n/2+3 \delta_{max}^2 \sqrt{3}/(32\, \omega)\le 0$, the average mean fitness $\langle \mb_\infty\rangle$ is a decaying function of $\mu$; otherwise, if $-n/2+3 \delta_{max}^2 \sqrt{3}/(32\, \omega)> 0$, $\langle \mb_\infty\rangle$ reaches a minimum for some value of $\mu$ in $(0,\omega/\sqrt{3})$, and then becomes concave and reaches a maximum for some larger value of $\mu$. This type of dependence, with the occurrence of an optimal mutation parameter (here, leading to a higher value of $\langle \mb_\infty\rangle$), has already been described for a fluctuating environment in \cite{CarNad20}, based on numerical simulations (see their figure~2; in their case, the mean population over one period is represented).

In \cite{FigMir18}, the same example (with $n=1$) $\delta(t)=\delta_{max} \, \sin(\omega\, t)$ was inspired by an experiment on the bacterial pathogen \emph{Serratia marcescens}. The method in \cite{FigMir18} is based on large time small mutation limit, and therefore can only give an equivalent of $\mb(t)$ at large times. Their work do not focus on the mean fitness, but on the mean trait $\xb(t)$ and variance $\ts^2(t)$ (the general theory deals with $n$ dimensions, but this particular example is in 1D):
$$\xb(t):=\int_\R x\, q(t,x) \, dx \hbox{ and }\ts^2(t):=\int_\R x^2\, q(t,x) \, dx-\xb(t)^2.$$They show that:
$$\xb(t)\approx \frac{\varepsilon \, \delta_{max}}{\omega}\sin(\omega \, t-\pi/2) \hbox{ and }\ts^2(t)\approx \varepsilon \, \sqrt{2}.$$In our framework, in the case $n=1$, the mean fitness is given by:
$$\mb(t)=-\frac{1}{2}\int_\R (x-\delta(t))^2 \, q(t,x) \, dx =-\frac{1}{2}\lp \ts^2(t)+\xb(t)^2\rp+\delta(t)\, \xb(t) - \frac{\delta(t)^2}{2}.$$Thus, their results can be used to compute an approached formula for $\mb(t):$
\begin{equation}\label{eq:FigMir18}
      \mb(t)\approx -\frac{\mu}{2}-\frac{1}{2}\lp \frac{\delta_{max}}{\omega} \rp^2(\omega \, \sin(\omega \, t)+ \mu \, \cos(\omega \, t))^2.
\end{equation}
We note that this is fully consistent with the large time asymptotics given by our formula \eqref{eq:mbar_sinus1}, in a small mutation regime ($\omega/(\omega^2+\mu^2)$ is approached by $1/\omega$) and with $n=1$.

Other forms of $\delta(t)$ could be considered as well, leading to more or less complex expressions for $\mb(t)$. For instance, with $\delta(t)=\delta_{max} \sin^2(\omega \, t)$, $\mb(t)$ is given by \eqref{eq:mbar_cor}, with
\begin{equation} \label{eq:Hdelta_sinus}
H_\delta(t)=\delta_{max}\left[ \frac{1}{2}-\frac{1}{8 \omega^2 + 2\mu^2}\lp 2\,  \omega \, \, \mu \, \sin(2 \omega t)\, \tanh(\mu t)+\mu^2 \cos(2  \omega t)+\frac{4\, \omega^2}{\cosh(\mu \, t)}\rp\right],
\end{equation}
and the average value of $\mb(t)$ over one period converges to
\begin{equation} \label{eq:mean_lag_sinus}
    \langle \mb_\infty\rangle:=\lim_{t\to +\infty}  \frac{\omega}{\pi}  \int_t^{t+\pi/\omega}\mb(s) \, ds=-\mu\, \frac{n}{2}-\frac{\delta_{max}^2 \, \omega^2  }{16 \omega^2 + 4\mu^2}.
\end{equation}
\Li Again, we observe that  $\langle \mb_\infty\rangle$ is a decreasing function of the frequency $\omega/\pi$. In light of formulas~\eqref{eq:mean_lag_sinusa} and~\eqref{eq:mean_lag_sinus}, a natural conjecture is that higher frequencies always tend to impede adaptation in periodically fluctuating environments. In the case of bounded domains, for equations of the form~\eqref{eq:main_form_n}, theorem~1.2 in \cite{CarNad20} shows that, for $r_{max}$ large enough (so that persistence occurs), the average value of the total population $\rho(t)$ over one period converges as $t\to+\infty$ towards the principal eigenvalue $\lambda_1(\omega)$ of the time-periodic operator $\phi\mapsto \partial_t \phi -\frac{\mu^2}{2}\Delta\,\phi-r(t,\x)\phi$, and $\rho(t)$ converges to a periodic function. Using the relationship \eqref{eq:rhot} (Appendix~A) between $\rho(t)$ and $\rb(t)$, integrating over one period and passing to the limit $t\to+\infty$, we obtain that $\langle \mb_\infty\rangle=-\lambda_1(\omega)-r_{max}$. Then, theorem~1.1 in \cite{LiuLou19} implies that $\lambda_1(\omega)$ is an increasing function of the frequency $\omega/\pi$. Finally, this shows that in bounded domains, $\langle \mb_\infty\rangle$ is indeed a decreasing function of the frequency. In our case, proving this result would require further investigation of the general formula~\eqref{eq:formule_gale_mbar} of Corollary~\ref{cor:formule_mbar}.  \Bk

\Li
\paragraph{Shifting and periodically fluctuating optimum.} Assume that $\delta(t)=\delta_1(t)+\delta_2(t)$, with $\delta_1(0)=\delta_2(0)=0$. Formula~\eqref{eq:mbar_cor} implies that:
$$\mb(t)=-\mu \, \frac{n}{2} \, \tanh(\mu \,t) -\frac{1}{2}\lp H_{\delta_1}(t)-\delta_1(t)+H_{\delta_2}(t)-\delta_2(t) \rp^2+  R_0'(t),$$as $H_\delta$ is linear with respect to $\delta$. Thus, we get:
$$\mb(t)=-\mu \, \frac{n}{2} \, \tanh(\mu \,t) +L_1(t)+L_2(t)-(H_{\delta_1}(t)-\delta_1(t))\,(H_{\delta_2}(t)-\delta_2(t))+  R_0'(t),$$with $$L_i(t):=-\frac{1}{2}\lp H_{\delta_i}(t)-\delta_i(t)\rp^2$$the lag associated with $\delta_i$. Let us now consider the specific example of a periodically fluctuating and shifting optimum. Combing the results of Propositions~\ref{prop:linear} and \ref{prop:sinus}, we obtain the following result.
\begin{proposition}\label{prop:sinus_shift}
Assume that $\delta(t)=c\,t+\delta_{max} \, \sin(\omega\, t)$
for some $c\in\R$, $\omega \in \R_+^*$ and $\delta_{max}>0$.
The mean fitness is given by
\begin{align}\label{eq:mbar_sin_shift}
      \mb(t)=&-\mu \, \frac{n}{2} \, \tanh(\mu \,t) \nonumber\\&-\frac{c^2}{2\, \mu^2} \tanh^2(\mu \, t) -\frac{1}{2}\lp\frac{\delta_{max} \,  \omega}{\omega^2 + \mu^2}\rp^2 ( \omega \, \sin( \omega \, t)+\mu \, \cos(
       \omega \, t) \tanh(\mu \, t))^2 \\&-\frac{c}{\mu}\tanh(\mu\,t)\,\lp\frac{\delta_{max} \,  \omega}{\omega^2 + \mu^2}\rp ( \omega \, \sin( \omega \, t)+\mu \, \cos(
       \omega \, t) \tanh(\mu \, t))\nonumber\\&+R_0'(t)\nonumber.
\end{align}
The average value of $\mb(t)$ over one period converges to
\begin{equation} \label{eq:mean_lag_sinus_shift}
    \langle \mb_\infty\rangle:=\lim_{t\to +\infty}  \frac{\omega}{2\,\pi}  \int_t^{t+2\,\pi/\omega}\mb(s) \, ds=-\mu\, \frac{n}{2}-\frac{c^2}{2\,\mu^2}-\frac{\delta_{max}^2 \, \omega^2  }{4 \omega^2 + 4\mu^2}.
\end{equation}
\end{proposition}
The effect of the two simultaneous changes (constant speed shift and periodic oscillations) is therefore not additive: we observe the emergence of an additional term that changes sign. However, in  average, over one period $ \langle \mb_\infty\rangle$ is the sum of the mutation load and of the lag load induced by each movement independently.

With the growth rate $\rb(t):=r_{max}+\mb(t)$ and assuming that the population size satisfies $\rho'(t)=\rb(t)\, \rho(t)$ or $\rho'(t)=\rho(t)(\rb(t)- \rho(t))$ (as in Appendix~A), we can compute the critical shifting speed for persistence (by persistence, we mean that $\rho(t)\not\to0$ as $t\to+\infty$). Proposition~\ref{prop:sinus_shift} shows that $\rb(t)$ converges to a periodic function, thus persistence occurs if and only if $\langle \rb_\infty\rangle:=r_{max}+\langle \mb_\infty\rangle\ge 0$ (the inequality is strict if $\rho'(t)=\rho(t)(\rb(t)- \rho(t))$). Applying formula~\eqref{eq:mean_lag_sinus_shift}, we obtain the following formula for the critical shifting speed: \begin{equation} \label{eq:maxspeed_sinshift}
c^*=\mu \, \sqrt{2 \, r_{max} - \mu\, n-\frac{\delta_{max}^2 \, \omega^2  }{2 \omega^2 + 2\mu^2}},
\end{equation}
provided that $r_{max}$ is large enough so that  $2 \, r_{max} - \mu\, n-\frac{\delta_{max}^2 \, \omega^2  }{2 \omega^2 + 2\mu^2}\ge0$ (otherwise persistence never occurs). We observe that the critical speed decreases with the frequency and with the amplitude of the periodic fluctuations.
\Bk

\section{Numerical computations \label{sec:numm}}

In this section, we check the validity of our results, (i)~compared to the initial integro-differential equation~\eqref{eq:main_IDE_moving_env}, and (ii)~compared to stochastic individual-based simulations of a standard model of genetic adaptation. With these two comparisons, we test the accuracy of the diffusion approximation (\eqref{eq:main_IDE_moving_env} vs \eqref{eq:main_PDE_moving_env}) and the effect of neglecting the stochastic aspects of mutation and selection. See Supplementary Material (S1) for some details on the numerical computation of the solution of~\eqref{eq:main_IDE_moving_env}.

\paragraph{Description of a Wright-Fisher individual-based model (IBM) with moving optimum.}
We assume a constant population size $N$. Under the assumptions of the Fisher's geometrical model, each individual $i=1,\ldots,N$  is characterized by a phenotype $\x_i\in\R^n$. Its relative Malthusian fitness at time $t$ (exponential growth rate) is given by \eqref{eq:m(t,x)}, i.e., $m_i=-\|\x_i-\Oc(t)\|^2/2$ and its corresponding Darwinian fitness is $\exp(m_i)$ (geometric growth rate, a discrete time counterpart of the Malthusian fitness). We assume non-overlapping generations of duration $\delta_t=1.$ Each generation, selection and genetic drift are jointly simulated by the multinomial sampling of $N$ individuals from the previous generation, each with weight given by their Darwinian fitnesses. Mutations are then simulated by randomly drawing, for each individual, a Poisson number of mutations, with rate $U$. We use a classic Gaussian FGM: each single mutation has a random phenotypic effect $d\x$ drawn into a multivariate Gaussian distribution: $d \x\sim \mathcal{N}(0,\lambda I_n)$, where $\lambda>0$ is the mutational variance at each trait, and $I_n$ is the identity matrix of size $n\times n$. Multiple mutations in a single individual have additive effects on phenotype.

\paragraph{Parameter values.} In all cases, we take $n=3$ and $\lambda=0.005$. Based on arguments in \cite{GilHamMarRoq19,MarRoq16}, the WSSM (diffusion) approximation should apply for $U\gtrsim U_c:= n^2\, \lambda/4$, see Supplementary Material (S2) for more details. We assume here
that  $U=10 \, U_c$  (recall that $\mu=\sqrt{U\, \lambda}$). Smaller values of $U$ are considered in Supplementary material (S2). In the individual-based model, we assume an initially clonal distribution of the phenotypes, at the optimum $0$; to be consistent with this assumption, we take $Q_0\equiv 0$.

\begin{figure}
\center
\subfigure[$\delta(t)=c\, t$]{\includegraphics[width=0.49\textwidth]{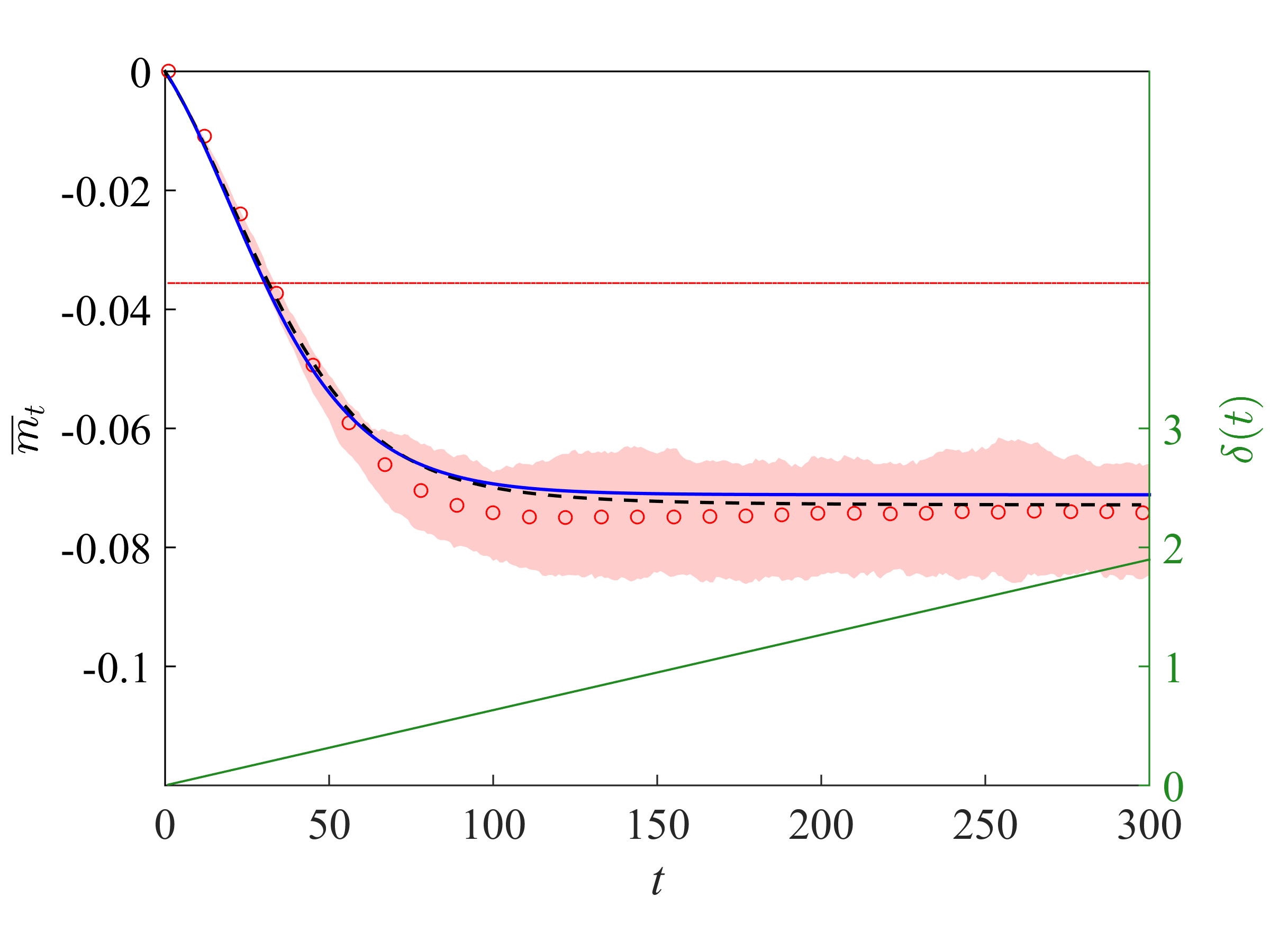}}
\subfigure[$\delta(t)=\delta_{max} \sin(\omega\, t)$]{\includegraphics[width=0.49\textwidth]{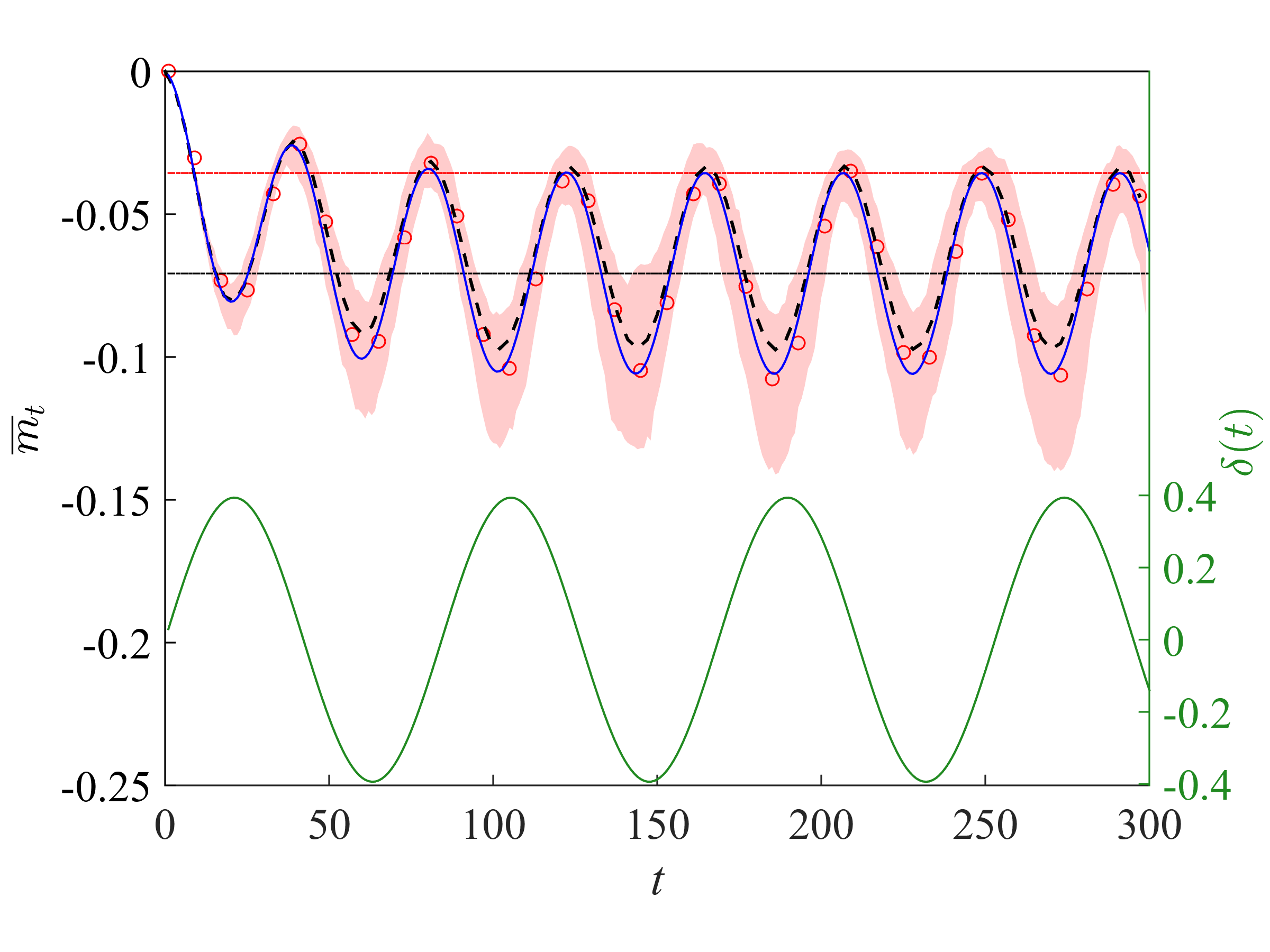}}
\subfigure[$\delta(t)=\delta_{max} \sin^2(\omega\, t)$]{\includegraphics[width=0.49\textwidth]{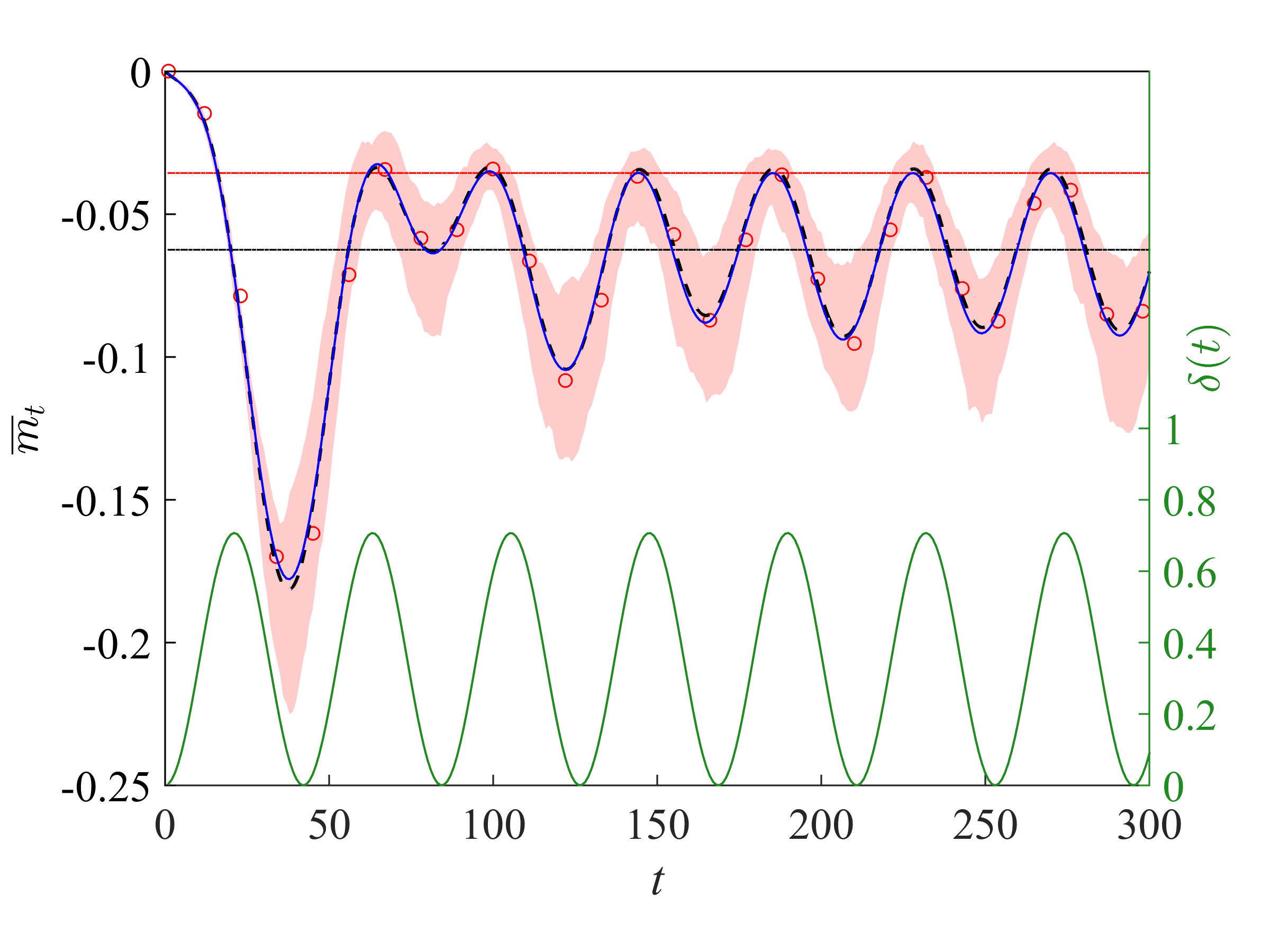}}
\subfigure[$\delta(t)=c\,t+\delta_{max} \sin(\omega\, t)$]{\includegraphics[width=0.49\textwidth]{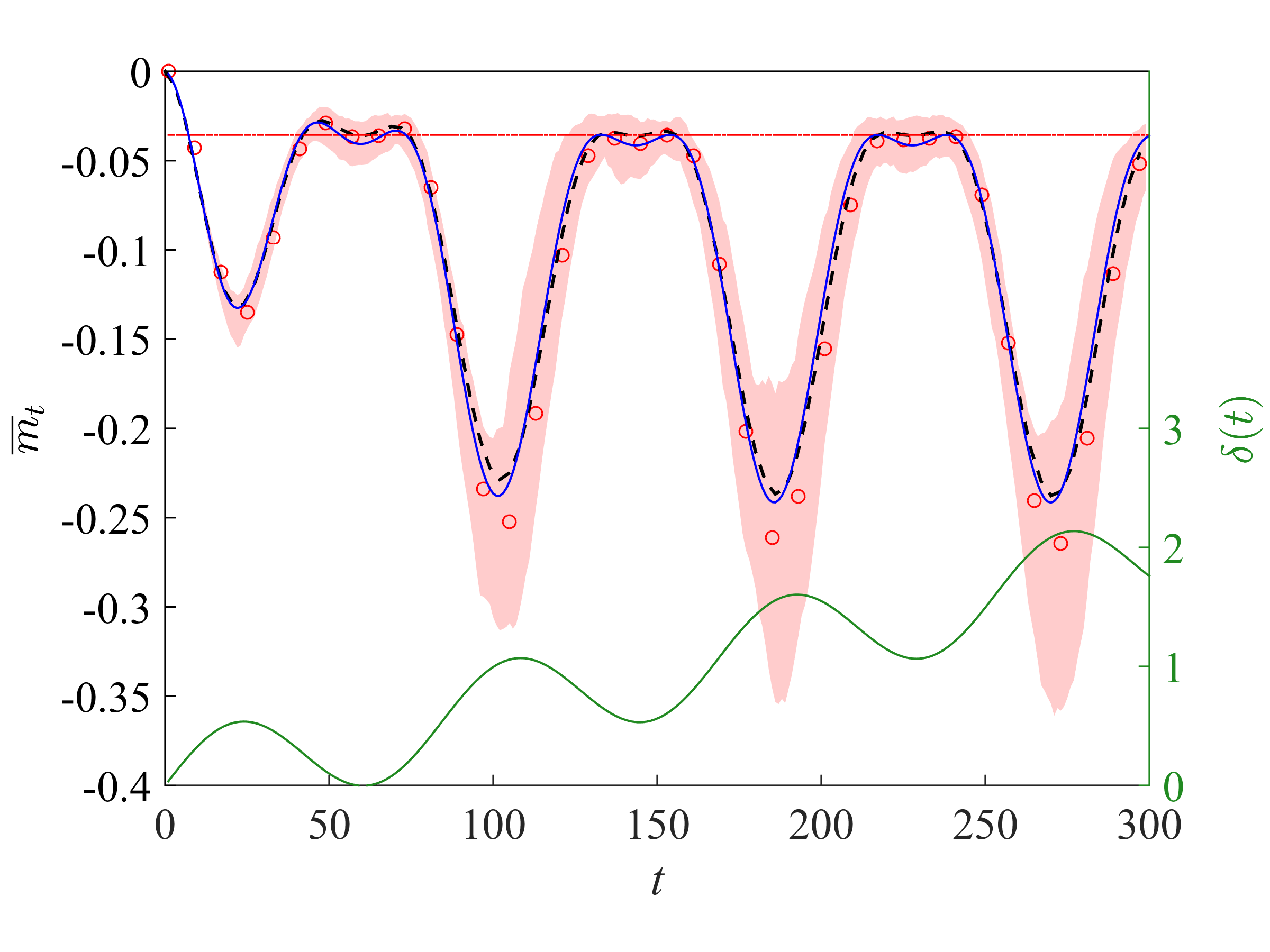}}
\caption{{\bf \Li Trajectories of mean fitness.\Bk}\Li \ Blue curve: theoretical value of $\mb(t)$ (given by Proposition~\ref{prop:linear} in panel a, by Proposition~\ref{prop:sinus} in panel b, by formula \eqref{eq:mbar_cor} with $H_\delta$ given by  \eqref{eq:Hdelta_sinus} in panel c and Proposition ~\ref{prop:sinus_shift} in panel d). Black dashed curve: numerical value obtained by solving the integro-differential equation~\eqref{eq:main_IDE_moving_env}; red circles: mean value of the mean fitness, averaged over $10^3$ replicate individual-based simulations ($N=10^4$ individuals in panel a; $N=10^3$ individuals in panels b,c,d); pink shading: interval between the 0.025 and 0.975 quantiles of the distribution of $\mb(t)$ obtained from the individual-based model; thin horizontal red
 line: mutation load $-\mu \, n /2$; thin horizontal black line (panels b,c): asymptotic average value of the mean fitness, $\langle \mb_\infty\rangle$ given by~\eqref{eq:mean_lag_sinusa} and~\eqref{eq:mean_lag_sinus};  green lines: $\delta(t)$. Parameter values: panel a: $c=\sqrt{n\, \mu^3}$; panel b:  $\delta_{max}=\sqrt{31\,\lambda}$ and $\omega= \mu \, \pi$; panel c:
 $\delta_{max}=10 \, \sqrt{\lambda}$ and $\omega= \mu \, \pi$; panel d: same values as in panels a,b.\Bk}
 \label{fig:linsin}
\end{figure}

\paragraph{Numerical results.} Fig.~\ref{fig:linsin}a depicts the trajectory of mean fitness when the optimum moves at a constant speed. We observe a good agreement between the analytical result of  Proposition~\ref{prop:linear}, the numerical value of $\mb(t)$ given by solving \eqref{eq:main_IDE_moving_env} and the mean value of $\mb(t)$ averaged over $10^3$ realisations of the IBM. Note that the value of the speed $c$ was chosen here such that the lag load is equal to the mutation load.

The trajectories of mean fitness corresponding to periodically varying optimums are presented in\Li \ Figs.~\ref{fig:linsin}b,c.  Again, the theoretical formulas accurately describe the average dynamics of the IBM and of the integro-differential equation~\eqref{eq:main_IDE_moving_env}. In particular, they capture the transient dynamics of adaptation, before $\mb(t)$ tends to become periodic. In Fig.~\ref{fig:linsin}c, the lowest value of $\mb(t)$ is reached during this transient stage, which means that extinction (or evolutionary rescue) will mainly depend on the early adaptation of the population, and not on the ultimate periodic behavior. In these plots, the parameter values are chosen such that the lag load averaged over one period is approximately equal to the mutation load.

In Fig.~\ref{fig:linsin}d, we considered the case of a shifting and periodically fluctuating optimum, corresponding to the situation studied in Proposition~\ref{prop:sinus_shift}. As expected, the trajectory of mean fitness is not just a combination of the trajectories of Figs.~\ref{fig:linsin}a,b: the extra term $$-\frac{c}{\mu}\tanh(\mu\,t)\,\lp\frac{\delta_{max} \,  \omega}{\omega^2 + \mu^2}\rp ( \omega \, \sin( \omega \, t)+\mu \, \cos(
       \omega \, t) \tanh(\mu \, t))$$in~ \eqref{eq:mbar_sin_shift} tends to lower some fluctuations and to increase others, leading to minimum values of $\mb(t)$ much lower than expected by simply adding the trajectories in panels a,b.
\Bk

\begin{figure}[h!]
\center
\subfigure[$U=U_c$]{\includegraphics[width=0.49\textwidth]{{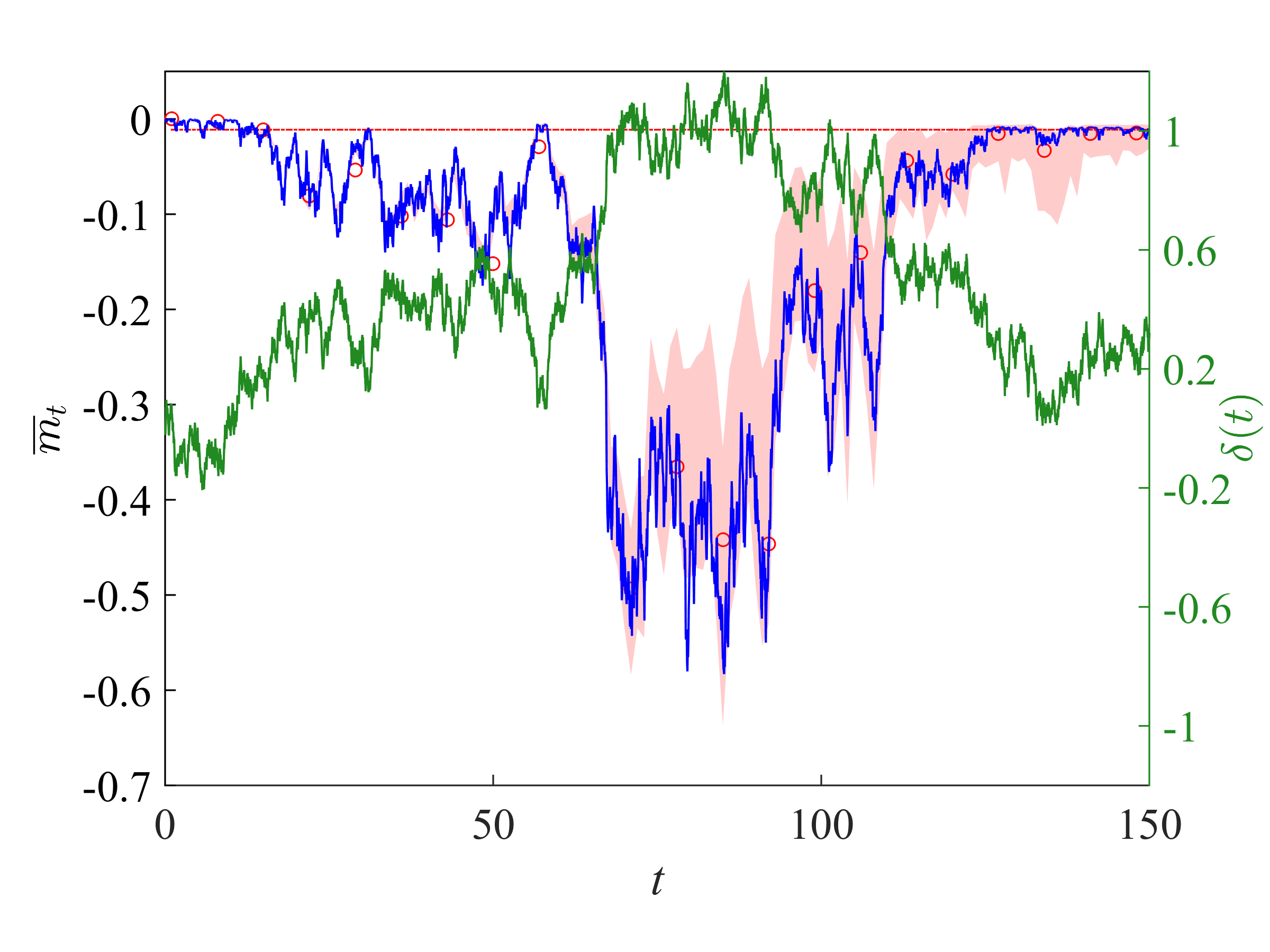}}}
\subfigure[$U=10\, U_c$]{\includegraphics[width=0.49\textwidth]{{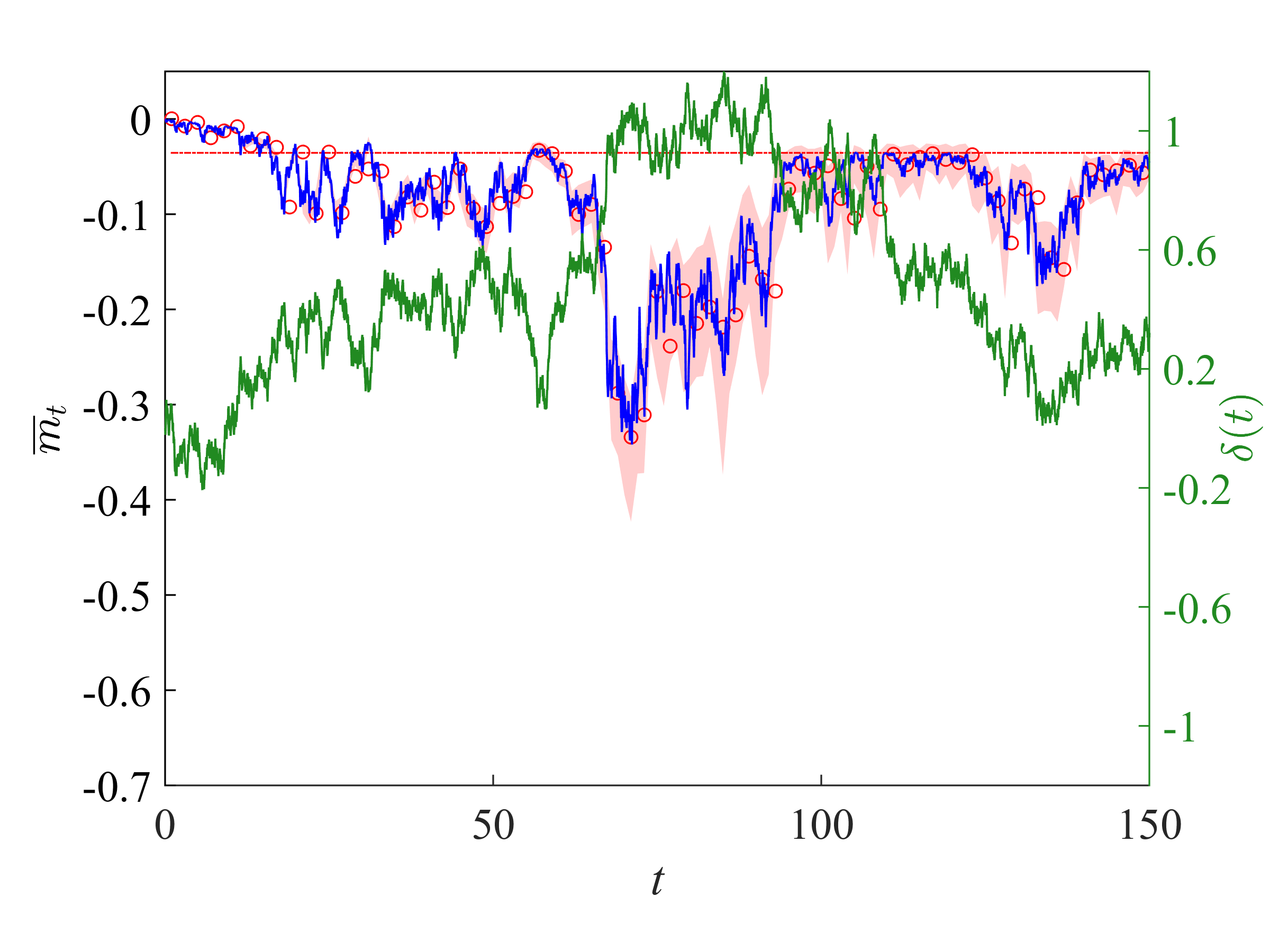}}}
\caption{{\bf Trajectories of mean fitness with a stochastically varying optimum.} Blue curve: theoretical value of $\mb(t)$ given by~\eqref{eq:mbar_cor}; red circles: mean value of the mean fitness, averaged over $10^3$ replicate individual-based simulations with $N=10^3$ individuals; pink shading: interval between the 0.025 and 0.975 quantiles of the distribution of $\mb(t)$ obtained from the individual-based model; thin horizontal red
 line: the mutation load $-\mu \, n /2$. We assumed here that the position of the optimum was given by $\Oc(t)=\delta(t)\, \bu$, with $\bu=(1,0,0)$ a unit vector in $\R^3$. Green line: $\delta(t)$, a realization of the Ornstein-Uhlenbeck process~\eqref{eq:OU_process}, with $\nu=0.01$ and $\tts=0.1$.}
 \label{fig:sto}
\end{figure}

Lastly, we tested the accuracy of the general formula~\eqref{eq:mbar_cor} in the case of a stochastic moving optimum. We assumed here that $\delta(t)$ was an Ornstein-Uhlenbeck process:
\begin{equation}\label{eq:OU_process}
d \delta(t)=- \nu \delta(t) \, dt + \tts \,  d W_t,
\end{equation}
with $W_t$ the Wiener process. Given a realization of this process, the formula~\eqref{eq:mbar_cor} can still be used to compute the value of $\mb(t)$ (though it requires a numerical evaluation of the integral in $H_\delta$). The results are presented in Fig.~\ref{fig:sto}. Again, the dynamics of the mean fitness simulated by the IBM are well-described by our theory. Note that all of the simulations were carried out based on a single realization $\delta(t)$ of the Ornstein-Uhlenbeck process. The comparison between Figs.~\ref{fig:sto} a) and b) illustrates the complex interplay between the environment and the mutation rate: the same environment leads to very different dynamics of adaptation depending on $U$ (or equivalently $\mu$).

\section{Discussion}

The approach developed in this paper led to an explicit characterization of the mean fitness and variance associated with the solution of \eqref{eq:main_PDE_moving_env} under very general assumptions on the movement of the optimum. Our results encompass in a single framework several examples that have been recently treated in the literature, and are consistent with these anterior results.

As pointed out in \cite{GomKro17}, most experimental studies are by nature  restricted to finite time horizons. Contrarily to `traveling wave', Hamilton-Jacobi,  `Gaussian solution' or spectral approaches, which deal with large time asymptotics, our framework leads to a description of the full dynamics of the mean fitness, since the initial state $t=0$, which is therefore particularly relevant for the understanding of drug resistance. Besides, as illustrated in Section~\ref{sec:numm}, the transient dynamics can be very different from the large time dynamics, even in a periodically fluctuating environment.

As mentioned in the Introduction the most widely used version of the FGM assumes that the Malthusian fitness decreases quadratically  away from the optimum~\cite{MarLen15,Ten14}. Several empirical tests showed that a deviation from the quadratic would actually lead to a lower fit of the data (see \cite{FraWel19} and fig. 5 in \cite{MarLen06a}). Thus, this assumption seems to be justified from a biological point of view. From a mathematical viewpoint, extending our approach to more general phenotype to fitness landscape model remains an interesting perspective. Being able to deal with anisotropic mutation effects, as in \cite{HamLav19} in the case of a fixed environment, could also lead to nontrivial results regarding the interaction between the direction of the environmental change and the type of anisotropy.

%Here, we briefly presented a preliminary approach to a full treatment of a stochastic environment, by showing the accuracy of our approach over a fixed realization of a stochastic process. An extension could be to derive the distribution of $\mb(t)$ in a full stochastic setting, i.e., when $\delta(t)$ is a random variable.

\section{ Proofs\label{sec:proofs}}

\subsection{Proof of Theorem~\ref{theo equation p}}
For any couple of functions $f$ and $g$ in $L^2(\R \times \R_-)$, we denote by $ \langle \cdot , \cdot  \rangle_{\R \times \R_-}$ the usual inner product on $L^2(\R \times \R_-)$:$$\langle f , g  \rangle_{\R \times \R_-}=\int_{\R \times \R_-} f(s_1,s_2) \, g(s_1,s_2) \, ds_1 \, ds_2.$$Similarly,  $ \langle \cdot , \cdot  \rangle_{\R^n}$ is the usual inner product on $L^2(\R^n)$.

Let $q(t,\x)$ be the solution of \eqref{eq:main_PDE_moving_env} defined in Theorem~\ref{thm:existuniq} and $\m(\x): \ \R^n \to \R \times \R_{-}$ the vector field defined by \eqref{components fitness}. For the sake of simplicity, for each $t\ge 0$, we denote by $q_t(\cdot)$ the function $\x\mapsto q(t,\x)$. We consider $\mathcal T$, the linear form defined over the Hilbert space $L^2(\R \times \R_-)$ by:$$
\mathcal T (\phi) :=  \langle  q_t, \phi \circ \m \rangle_{\R^n}.$$First, we ensure that $ \mathcal T $ is a well-defined continuous linear form. Let $\psi \in L^1(\R^n)$. Thanks to  a change of variables, the following identity holds true:
\begin{multline}\label{change variable}
\int_{\R^n} \psi\left( x_1,- \frac{\| \x \|^2}{2},x_3,\ldots ,x_n\right) d\x =  \\
2 \int_{E_{1,2}}  \left( \int_{E_{n-2}} \psi(y_1,y_2,\ldots ,y_n) \frac{1}{\sqrt{-2y_2 -y_1^2 -  \sum_{i=3}^n y_i^2 }}dy_3 \ldots dy_n \right) dy_1dy_2,
\end{multline}
where:
\begin{align*}
& E_{1,2} = \Big\{ (y_1,y_2) \in \R \times \R_- / \, -2y_2 -y_1^2\geq 0  \Big\}, \text{ and }\\
& E_{n-2} = \Big\{ (y_3,\ldots ,y_n) \in \R^{n-2} / \,   \sum_{i=3}^n y_i^2 \leq    -2y_2 -y_1^2  \Big\} \text{ for } (y_1,y_2) \in E_{1,2}.
\end{align*}
We will use this formula to show that $q_t(\cdot) \phi \circ \m(\cdot)$ is a $L^1(\R^n)$ function and that $\mathcal T$ is a continuous linear form. Let $\phi \in L^2(\R \times \R_-)$. The Cauchy-Schwartz inequality implies that, for any $b >0$,:
\begin{align}\label{eq:CS T}
\abs{ \mathcal T(\phi) } & = \abs{ \int_{\R^n} q_t(\x) \phi(\m(\x)) d\x  } \nonumber \\ & \leq   \abs{\int_{\R^n} \phi \left( x_1,- \frac{\| \x \|^2}{2}\right)^2 e^{-b \| \x \|} d\x}^{\frac 12} \abs{\int_{\R^n}q_t(\x)^2 e^{b \| \x \|} d\x}^{\frac 12}.
\end{align}
Using \eqref{eq:q_exp_bound} we know that $q_t$ is exponentially bounded, and so is $q_t^2$. Thus the last term of the above inequality is finite. We can now deal with the other part. Applying \eqref{change variable} with $\psi(\x)=\phi(x_1,-\|\x\|^2/2)^2\, \exp(-b \, \|\x\|)$, we get
\begin{multline}\label{eq:int cont T}
\int_{\R^n} \phi \left( x_1,- \frac{\| \x\|^2}{2}\right)^2 e^{-b \| \x\|} d\x=  \\
2 \int_{E_{1,2}} \phi(m_1,m_2)^2 e^{-b \sqrt{-2 m_2}}  \left( \int_{E_{n-2}}  \frac{1}{\sqrt{-2m_2 -m_1^2 -  \sum_{i=3}^n y_i^2 }}dy_3 \ldots dy_n \right) dm_1dm_2.
\end{multline}
Moreover, by a polar change of coordinates, we know that, for any $y\geq0$, and $k \in \N$:
\begin{align*}
\int_{ \sum_{i=1}^k y_i^2 \leq y } \frac{1}{\sqrt{y -  \sum_{i=1}^k y_i^2 }}dy_1 \ldots dy_k = \int_0^{ \sqrt{y}} \frac{1}{\sqrt{y - r^2 }} \hbox{Vol}_{k-1}(\sqrt{y}) dr,
\end{align*}
where $\hbox{Vol}_{k-1}(s)$ is the volume of the ball of radius $s$ in dimension $k-1$, that is $V_k(s)= C_{k-1} s^{k-1}$, with $C_{k-1}$ a generic constant depending only upon the dimension $k$. Therefore,
\begin{align*}
\int_{ \sum_{i=1}^k y_i^2 \leq y } \frac{1}{\sqrt{y -  \sum_{i=1}^k y_i^2 }}dy_1 \ldots dy_k & = \frac{\pi}{2} \hbox{Vol}_{k-1}(\sqrt{y}), \\
& = C_k y^{\frac {k-1}2} .
\end{align*}
Plugging this computation into \eqref{eq:int cont T}, we have found
\begin{align*}
\int_{\R^n} \phi \left( x_1,- \frac{\| \x \|^2}{2}\right)^2 e^{-b \| \x\|} d\x & =    C_{n-2} \int_{E_{1,2} } \phi(m_1,m_2)^2e^{-b \sqrt{-2 m_2}}   (-2 m_2 -m_1^2)^{\frac{n-3}{2}}  dm_1dm_2 ,\\
& \leq C_n \int_{E_{1,2} } \phi(m_1,m_2)^2dm_1dm_2  \\&\leq C_n \int_{\R \times \R_- } \phi(m_1,m_2)^2dm_1dm_2,
\end{align*}
for some generic constant $C_n$ that depends only on the dimension $n$. Using this estimate together with \eqref{eq:CS T}, we get:
\begin{align*}
\abs{\mathcal T (\phi)}\leq C_q \norm{\phi}_{\Li L^2 \Bk(\R\times \R_-)},
\end{align*}
where $C_q$ depends only on the bound on $q$ established in \eqref{eq:q_exp_bound} and on the dimension. Therefore, $\mathcal{T}$ is a well-defined \emph{continuous} linear form on $\Li L^2 \Bk(\R\times \R_-)$.

From then, the existence and uniqueness of
a function $p_t\in L^2(\R\times \R_-)$ satisfying $$\mathcal{T}(\phi)=\langle p_t , \phi \rangle_{\R\times \R_-}$$for all $\phi \in L^2(\R \times \R_-)$
is  a straightforward application of the Riesz-Frechet representation \Li theorem, \Bk see for instance theorem 6.19 in~\cite{Rud06}. Defining $p(t,\cdot,\cdot):=p_t$ for each $t\ge 0$, the regularity $p\in C^1(\R_+,L^2(\R \times \R_-))$ is a straightforward consequence of the regularity of $q$. $\Box$

\subsection{Proof of Theorem~\ref{thm:CGF}}
Fix $z_1,z_2 \in \R \times \R_+$ and consider an \emph{increasing} sequence $(\phi_k)_{k\in \N}$ of nonnegative functions in $\mathcal{C}_c^{\infty}(\R \times \R_-)$ satisfying:
\begin{equation} \label{def:phik}
    \left\{\begin{array}{l}
         \phi_k(m_1,m_2)=\exp(z_1 m_1 + z_2 m_2),  \hbox{ if }m_1^2+m_2^2<k,\\
         \phi_k(m_1,m_2)=0,  \hbox{ if }m_1^2+m_2^2>k+1.
    \end{array}\right.
\end{equation}
Applying Theorem~\ref{theo equation p},
one gets
\begin{align*}
\int_{\R^n} q(t,\x) \phi_k(\m(\x)) d\x = \int_{\R\times \R_-} p(t,m_1,m_2) \phi_k(m_1,m_2)dm_1dm_2,
\end{align*}
with $\m(\x)=(\mm_1(\x),\mm_2(\x)): \ \R^n\to \R\times \R_-$ defined by \eqref{components fitness}. Using the monotone convergence theorem, we can pass to the limit $k\to +\infty$ on each side of the above equality. This yields:
\begin{align*}
\int_{\R^n} q(t,\x)  \exp \left( z_1 \bu \cdot \x- z_2 \frac{\norm{\x}^2}{2}  \right) d\x = \int_{\R\times \R_-} p(t,m_1,m_2) \exp (z_1 m_1 +z_2 m_2)dm_1dm_2.
\end{align*}
According to \eqref{eq:q_exp_bound}, $q$ is integrable against exponential functions. Therefore, the left hand side is finite (and positive, since $q\not \equiv 0$), and so is the right hand side. Thus, the quantity $$C(t,z_1,z_2):=\ln \lp \int_{\R \times \R_{-}}p(t,m_1,m_2)\, e^{m_1 \, z_1 +m_2 \, z_2} \, dm_1 \, dm_2 \rp,$$is well-defined.

\Rd
We now detail how to derive the equation \eqref{eq:C} satisfied by $C$. First,  differentiating \eqref{link q and p} with respect to time, we find that for any compactly supported test function $\phi\in \C^\infty_c(\R \times \R_-,\R)$:
\begin{align*}
\langle \p_t p, \phi   \rangle_{\R \times \R_{-}} = \langle  \p_t q , \phi \circ \m \rangle_{\R^n}.
\end{align*}
Plugging in \eqref{eq:main_PDE_moving_env}, the equation solved by $q$, for any $t >0$, one has
\begin{align}\label{eq:ptp}
\langle \p_t p(t, \cdot,\cdot) , \phi(\cdot,\cdot)  \rangle_{\R \times \R_{-}} = \frac{\mu^2}{2}\langle \Delta q(t,\cdot) , \phi ( \m(\cdot)) \rangle_{\R^n} + \langle q(t,\cdot) (m(t,\cdot)- \mb(t)) , \phi ( \m( \cdot) )\rangle_{\R^n} .
\end{align}
We deal separately with each term on the right hand side of \eqref{eq:ptp}.
First, for the Laplace operator, we use a duality argument:
\begin{align*}
 \frac{\mu^2}{2}\langle \Delta q  , \phi \circ \m \rangle_{\R^n} =  \frac{\mu^2}{2}\langle  q ,  \Delta( \phi \circ  \m )  \rangle_{\R^n} ,
\end{align*}
since $\Delta$ is self-adjoint in $L^2(\R^n)$. We can then write:
\begin{align}\label{eqint_Last}
\frac{\mu^2}{2}\Delta (\phi \circ \m)  = (\mathcal A \phi) \circ \m,
\end{align}
with $\mathcal{A}$ such that for all  $(m_1,m_2)$ in $\R \times \R_{-},$
\begin{multline}\label{explicit Abis}
\mathcal A (\phi)(m_1,m_2) =  -\mu^2 m_2 \p_{22} \phi(m_1,m_2) + \frac{\mu^2}{2}  \p_{11} \phi(m_1,m_2)  -\mu^2 m_1 \p_{21} \phi(m_1,m_2)  \\
-\frac{\mu^2 n}{2} \p_{2} \phi(m_1,m_2).
\end{multline}
With $p$ defined by Theorem~\ref{theo equation p}, we have:
\begin{align}
\nonumber \frac{\mu^2}{2} \langle \Delta q  , \phi \circ \m \rangle_{\R^n} & = \langle  q ,   (\mathcal A \phi) \circ \m \rangle_{\R^n} ,\\
& \label{eq:ptp int0} = \langle  p ,   \mathcal A (\phi) \rangle_{\R \times \R_-}.
\end{align}
Going back to \eqref{eq:ptp}, we notice that thanks to  \eqref{eq:m(t,x)_b},
\begin{align*}
\langle q(t,\cdot) m(t,\cdot) , \phi ( \m( \cdot) ) \rangle_{\R^n}  = \left\langle q(t,\cdot)  \left( \delta(t)\, \mm_1(\cdot)+\mm_2(\cdot)-\frac{\delta(t)^2}{2} \right) , \phi ( \m( \cdot) ) \right\rangle_{\R^n}  , \end{align*}
which can be decomposed as
\begin{multline*}
 \langle q(t,\cdot) m(t,\cdot) , \phi ( \m( \cdot) ) \rangle_{\R^n}  =  \delta(t)  \left\langle q(t,\cdot)  ,\mm_1(\cdot)\,\phi ( \m( \cdot) ) \right\rangle_{\R^n}   + \left\langle q(t,\cdot)  , \mm_2(\cdot)\,\phi(\m(\cdot)) \right\rangle_{\R^n}  \\ -\frac{\delta(t)^2}{2} \left\langle q(t,\cdot)  , \phi ( \m( \cdot) ) \right\rangle_{\R^n}.
\end{multline*}
We apply Theorem \eqref{theo equation p} three times, to get
\begin{align}\label{eq:ptpt int}
\langle q(t,\cdot) m(t,\cdot) , \phi ( \m( \cdot) ) \rangle_{\R^n} & = \left\langle p(t,m_1,m_2) ,   \left( \delta(t)\, m_1 +m_2 -\frac{\delta(t)^2}{2} \right)  \phi ( m_1, m_2 ) \right\rangle_{\R \times \R_- } .
\end{align}
To deal with the last part of \eqref{eq:ptp}, we first notice that:
\begin{align*}
\mb(t)=\langle q(t,\cdot)  , m(t, \cdot) \rangle_{\R^n} &= \left\langle q(t,\cdot)  , \delta(t)\, \mm_1(\cdot)+\mm_2(\cdot)-\frac{\delta(t)^2}{2}  \right\rangle_{\R^n}\\ & =\delta(t) \langle p(t,m_1,m_2) ,m_1 \rangle_{\R \times \R_- }   + \langle p(t,m_1,m_2) ,m_2 \rangle_{\R \times \R_- } -\frac{\delta(t)^2}{2}.
\end{align*}
Again, the above equality involves Theorem \eqref{theo equation p}, with $\phi=\hbox{Id}_{\R \times \R_-}$, which is made possible by considering increasing compactly supported approximations of identity, and a uniform limit as in \eqref{def:phik}. Thus,
\begin{align}\label{eq:ptpt int2}
\mb(t) \langle q(t,\cdot)  , \phi ( \m( \cdot) ) \rangle_{\R^n} & = \left( \delta(t)\, \hm_1(t)+ \hm_2(t)-\frac{\delta(t)^2}{2} \right) \left\langle p(t,m_1, m_2 )   , \phi ( m_1,m_2) ) \right\rangle_{\R \times \R_- }  ,
\end{align}
where $\hm_j$ is defined by $\hm_j(t) :=\langle p(t,m_1,m_2) ,m_j \rangle_{\R \times \R_- }$  ($j=1,2$).
Plugging \eqref{eq:ptp int0}, \eqref{eq:ptpt int} and \eqref{eq:ptpt int2} into \eqref{eq:ptp}, we find that
\begin{multline}\label{eq pweak}
\langle \p_t p(t, m_1,m_2) , \phi(m_1,m_2)  \rangle_{\R \times \R_{-}} = \left\langle p (t, m_1, m_2) , \mathcal A (\phi)(m_1,m_2)  \right\rangle_{\R \times \R_{-}}  + \\
\left\langle p(t,m_1,m_2) ,   \Big( \delta(t)\, (m_1 -\hm_1(t) ) +m_2- \hm_2(t)  \Big)   \phi ( m_1, m_2 ) \right\rangle_{\R \times \R_- }.
\end{multline}
To be able to compute an equation on the cumulant generating function $C$ defined by \eqref{def:C(t,z1,z2)},  we first define the moment generating function, for all $z_1\in\R$ and $z_2\in\R_+$:
\begin{align*}
M(t,z_1,z_2) := \langle p ,  E_{z_1,z_2} \rangle_{\R \times \R_-} \text{ with } \,
 E_{z_1,z_2}(m_1,m_2) & := \exp(z_1 m_1 + z_2 m_2).
\end{align*}
Applying the relationship \eqref{eq pweak} to the sequence $(\phi_k)_{k\in \N}$ defined by \eqref{def:phik} and passing to the limit $k\to +\infty$, we get:
\begin{multline}
\langle \p_t p(t, \cdot,\cdot), E_{z_1,z_2}(\cdot,\cdot)   \rangle_{\R \times  \R_{-}} = \Big\langle   p(t,m_1,m_2) , \, \mathcal{A}(E_{z_1,z_2} )(m_1,m_2)  \\
+ \Big( \delta(t)\, (m_1 -\hm_1(t) ) +m_2- \hm_2(t)  \Big)    E_{z_1,z_2}(m_1,m_2) \Big\rangle_{\R \times \R_{-}}.
\end{multline}
With the explicit expression of $ \mathcal{A}$ given in \eqref{explicit Abis}, we can compute $\p_t M$ since
$$\p_t M(t,z_1,z_2) = \langle \p_t p(t, \cdot,\cdot), E_{z_1,z_2}(\cdot,\cdot)   \rangle_{\R \times  \R_{-}}.$$
We find,
\begin{multline}\label{eq:int calc C}
\p_t M(t,z_1,z_2) = \left\langle   p, -\mu^2 m_2 z_2^2 E_{z_1,z_2} +  \frac{\mu^2}{2} z_1^2 E_{z_1,z_2} -\mu^2 m_1 z_2z_1 E_{z_1,z_2} -\frac{\mu^2 n}{2} z_2 E_{z_1,z_2}\right\rangle_{\R \times \R_{-}} \\ + \left\langle   p,  E_{z_1,z_2}  \Big( m_2-\hm_2 + \delta(t)(m_1-\hm_1)  \Big) \right\rangle_{\R \times \R_{-}}.
\end{multline}
Moreover, for $j=1,2$,
 \begin{align*}
 \langle   p , m_j  E_{z_1,z_2} \rangle_{\R \times \R_{-}} = \p_j M(t,z_1,z_2) \hbox{ and } \hm_j(t)=\langle   p , m_j   \rangle_{\R \times \R_{-}} = \p_j M(t,0,0).
 \end{align*}
 Therefore, \eqref{eq:int calc C} can be rewritten as
 \begin{align*}
 \p_t M(t,z_1,z_2) =&   \delta(t) \Big(\p_{1} M(t,z_1,z_2)- M(t,z_1,z_2) \p_{1} M(t,0,0)  \Big)\\ &+
  \p_{2} M(t,z_1,z_2)-  M(t,z_1,z_2) \p_{2} M(t,0,0) -\mu^2 z_1 \, z_2 \p_{1} M(t,z_1,z_2) \\ &-\mu^2 z_2^2 \p_{2} M(t,z_1,z_2) + \mu^2 ( z_1^2/2 -n\, z_2 /2) M(t,z_1,z_2).
 \end{align*}
Dividing this expression by $M(t,z_1,z_2)$, and  since $C(t,z_1,z_2)= \ln(M(t,z_1,z_2))$, this shows that $C$ satisfies the equation in Theorem~\ref{thm:CGF}. $\Box$\Bk

\subsection{Proof of Theorem~\ref{th:formule_C}}

\Li We begin with the proof of Proposition~\ref{prop:solQ}.\Bk

\noindent \textit{Proof.}  Let $T>0$. For all $t\in [0,T]$ and $z, \, \tz \in [t-T,+\infty)^2$ we define $$W(t,z,\tz)=Q(t,z+T-t,\tz+T-t)-Q(t,T-t,T-t).$$First, we observe that:
\begin{multline}
\partial_t W(t,z,\tz)= (\partial_t Q - \p_z Q -\p_{\tz}Q) (t,z+T-t,\tz+T-t)\\ - (\partial_t Q - \p_z Q-\p_{\tz}Q) (t,T-t,T-t).
\end{multline}
Using \eqref{eq:Q_inia}, we then obtain:
\begin{equation*}
\partial_t W(t,z,\tz)= \beta(t,z+T-t,\tz+T-t)-\beta(t,T-t,T-t).
\end{equation*}
Integrating between $0$ and $t$, and noting that $W(0,z,\tz)=Q_0(z+T,\tz+T)-Q_0(T,T),$ we get:
$$
W(t,z,\tz)=\int_0^t \beta(s,z+T-s,\tz+T-s)-\beta(s,T-s,T-s)\, ds+W(0,z,\tz),
$$
which leads to:
\begin{multline}
 \label{eq:Q_int1} Q(t,z+T-t,\tz+T-t)-Q(t,T-t,T-t)\\=\int_0^t \beta(s,z+T-s,\tz+T-s)-\beta(s,T-s,T-s)\, ds+ W(0,z,\tz).\end{multline}
Computing this quantity at $(z,\tz)=(t-T,t-T)$ and using $Q(t,0,0)=0$, we get:
\begin{equation} \label{eq:Q_int2} -Q(t,T-t,T-t)=\int_0^t \beta(s,t-s,t-s)-\beta(s,T-s,T-s)\, ds+ W(0,t-T,t-T).\end{equation} Combining \eqref{eq:Q_int1} and \eqref{eq:Q_int2}, we obtain:
\begin{multline} \label{eq:Q_int3} Q(t,z+T-t,\tz+T-t)= \int_0^t  \beta(s,z+T-s,\tz+T-s)-\beta(s,t-s,t-s) \,ds\\ + W(0,z,\tz)- W(0,t-T,t-T),\end{multline}
which implies that $$Q(t,z,\tz)=
\int_0^t  \beta(s,z+t-s,\tz+t-s)-\beta(s,t-s,t-s) \,ds + Q_0(z+t,z+t)- Q_0(t,t).$$Conversely, it is straightforward to check that this expression solves \eqref{eq:Q_inia}, and the proposition follows. $\Box$

\

In order to solve our main equation \eqref{eq:C}, we look for a function $\phi_t(z,\tz)=(y_1(t,z,\tz),y_2(z))$,
with $y_1\in \C^1(\R_+^3,\R)$ and $y_2 \in \C^1(\R_+,\R_+)$,
such that the function
\begin{equation} \label{def:Q}
Q(t,z,\tz)=C(t,\phi_t(z,\tz))
\end{equation}
satisfies a problem of the form \eqref{eq:Q_inia}. In that respect, we first establish some conditions on the functions $y_1,$ $y_2$.

\begin{lem}\label{lem:sys_y}
Assume that $C$ is a solution of \eqref{eq:C}, and assume that $y_1\in \C^1(\R_+^3,\R)$ and $y_2 \in \C^1(\R_+,\R_+)$, satisfy, for $(t,z,\tz) \in \R_+^3$:
\begin{equation}
\label{eq:sys_y}
\left\{\begin{array}{l}
\partial_t y_1 -\partial_z y_1-\partial_{\tz} y_1=-\delta(t)+\mu^2 \,y_1\,y_2,  \\
\partial_z y_2=1- \mu^2\, y_2^2, \\
y_1(t,0,0)=0, \ \partial_t y_1(t,0,0)= 0,  \ y_2(0)=0,
\end{array} \right.
\end{equation}
 then the function $Q(t,z,\tz)$ defined by \eqref{def:Q} satisfies \eqref{eq:Q_inia}, with $\beta(t,z,\tz)=\gamma(y_1(t,z,\tz),y_2(z))$ and $Q_0(z,\tz)=C_0(y_1(0,z,\tz),y_2(z))$.
\end{lem}
\noindent\emph{Proof of Lemma~\ref{lem:sys_y}.} Let $Q$ be defined by \eqref{def:Q}. Then, one can note that:
\begin{equation}
\label{eq:deriv_Q1}
\left\{\begin{array}{l}
\partial_t Q(t,z,\tz)= \partial_t C(t,\phi_t(z,\tz))+ \partial_t y_1(t,z,\tz) \partial_1 C(t,\phi_t(z,\tz)), \\
\partial_z Q(t,z,\tz)= \partial_z y_1(t,z,\tz) \partial_1 C(t,\phi_t(z,\tz))+\partial_z y_2(z) \partial_2 C(t,\phi_t(z,\tz)),\\
\partial_{\tz} Q(t,z,\tz)= \partial_{\tz} y_1(t,z,\tz) \partial_1 C(t,\phi_t(z,\tz)).
\end{array} \right.
\end{equation}
Thus, $Q$ satisfies \eqref{eq:Q_inia}, with $\beta(t,z,\tz)=\gamma(\phi_t(z,\tz))$ if and only if
\begin{equation} \label{eq:eq_Cy}
\begin{array}{rl}
\partial_t C (t,\phi_t(z,\tz))  = &(\partial_z y_1+\partial_{\tz} y_1-\partial_t y_1) \partial_1 C(t,\phi_t(z,\tz)) + (\partial_z y_2) \, \partial_2 C(t,\phi_t(z,\tz)) \\
& -(\partial_z y_1+\partial_{\tz} y_1)(t,0,0)\, \partial_1 C(t,\phi_t(0,0))-\partial_z y_2(0)\, \partial_2 C(t,\phi_t(0,0)) \\
& +\gamma(\phi_t(z,\tz)).\end{array}
\end{equation}
Using \eqref{eq:sys_y}, we just have to check that the coefficients in front of the differential terms in \eqref{eq:eq_Cy} correspond to those in \eqref{eq:C}, computed at $(t,\phi_t(z,\tz))$ to conclude the proof.
Note that, at $z=0$, $\partial_t y_1 -\partial_z y_1-\partial_{\tz} y_1=-\delta(t)$, and  $\partial_t y_1(t,0,0)= 0$ thus implies that $(\partial_z y_1+\partial_{\tz} y_1)(t,0,0)=\delta(t)$. $\Box$

\

Our goal is now to find some functions $y_1$, $y_2$, satisfying the conditions of Lemma~\ref{lem:sys_y}. In order to solve the system \eqref{eq:sys_y}, we first note that
\begin{equation}\label{eq:soly2}
y_2(z)=\tanh(\mu \, z)/\mu
\end{equation}
satisfies the second equation in the system and the condition $y_2(0)=0$. Then, fix $T>0$ and define, for $t\in[0,T)$ and $z,\, \tz \in [t-T,+\infty)^2$, $h(t,z,\tz):=y_1(t,z+T-t,z+\tz+T-t).$
The function $h$ satisfies
\begin{equation}\label{eq:h}
\partial_t h(t,z,\tz) =-\delta(t)+\mu^2\, h(t,z,\tz) \, y_2(z+T-t),
\end{equation}
which can be solved explicitly. Namely, for any  function $B$ in $C^1(\R^2)$, a solution is given by:$$h(t,z,\Li\tz\Bk)=F(t,z+T-t)\lp B(z+T,\tz)- \int_0^t \delta(s)F(-s,z+T)\, ds \rp,$$with
$$F(t,z)=\frac{\cosh(\mu(z+t))}{\cosh(\mu \, z)}.$$This leads to
 the following expression for $y_1(t,z,\tz)=h(t,z+t-T,\tz-z)$:
\begin{equation} \label{eq:solgaley1}
y_1(t,z,\tz)=F(t,z)\lp B(z+t,\tz-z)- \int_0^t \delta(s)F(-s,z+t)\, ds  \rp.
\end{equation}
The function $B$ must be such that $y_1(t,0,0)=0$. We chose
\begin{equation*}
B(z+t,\tz-z)=z-\tz+\int_0^{z+t} \delta(s)F(-s,z+t)\, ds,
\end{equation*}
and finally get:
\begin{equation} \label{eq:solgaley1b}
y_1(t,z,\tz)= \int_0^{z} \delta(z+t-s) \, \frac{\cosh(\mu \, s)}{\cosh(\mu \, z)}\, ds+(z-\tz)\, F(t,z).
\end{equation}
Finally, it is immediate to check that $y_2$ and $y_1$ respectively defined by \eqref{eq:soly2} and \eqref{eq:solgaley1b} satisfy the conditions of Lemma~\ref{lem:sys_y}.

Now, let $C\in \C^{1} (\R_+\times\R \times \R+,\R)$ be defined by \eqref{eq:C}, $\phi_t(z,\tz):=(y_1(t,z,\tz),y_2(z)): \ \R_+ \to \R\times \R_+$ and set $Q(t,z,\tz)=C(t,\phi_t(z,\tz))$ for all $t\ge 0$ and $z,\, \tz \in \R_+^2$. Then, Lemma~\ref{lem:sys_y} implies that $Q(t,z,\tz)$ satisfies~\eqref{eq:Q_inia} with $Q_0(z,\tz)=C_0(\phi_0(z,\tz))$ and $\beta(t,z,\tz)=\gamma(\phi_t(z,\tz))$. Proposition~\ref{prop:solQ} implies that $Q(t,z,\tz)$ is given by the expression~\eqref{eq:solQ}. This proves the result of Theorem~\ref{th:formule_C}. $\Box$

\subsection{Proof of Corollary~\ref{cor:formule_mbar}}
Set $R(t,z)=Q(t,z,z)$. We have $\p_z R(t,0)=\p_z Q(t,0,0)+\p_{\tz} Q(t,0,0)$. Using \eqref{eq:deriv_Q1}, we observe that $$\partial_z R(t,0)=(\p_z y_1+\p_{\tz} y_1)(t,0,0) \partial_1 C(t,0,0)+\p_z y_2(0) \p_2 C(t,0,0).$$From \eqref{eq:sys_y}, we have $(\p_z y_1+\p_{\tz} y_1)(t,0,0)=\delta(t)$ and $\p_z y_2(0)=1$. Finally, $\partial_z R(t,0)=\delta(t)\, \p_1 C(t,0,0)+\p_2 C(t,0,0)$, and using \eqref{eq:mbar_CGF} we obtain the general formula \eqref{eq:formule_gale_mbar} for $\mb(t):$
\begin{equation} \label{eq:mb_gale1}
\mb(t)=\partial_z R(t,0)-\frac{\delta(t)^2}{2}.
\end{equation}
To derive a more explicit expression, we begin by observing that:
\begin{equation} \label{eq:mb_beta}
\partial_z R(t,0)=\int_0^t  (\partial_z \beta+\partial_{\tz} \beta)(t-u,u,u) \, du + R_0'(t),
\end{equation}
and
\begin{equation} \label{eq:beta}
\begin{array}{rl}
      (\partial_z \beta+\partial_{\tz} \beta)(t-u,u,u)  & = \ds  \left[(\partial_z y_1 +\partial_{\tz} y_1) \p_1 \gamma(y_1,y_2) + \p_z \, y_2 \, \p_2 \gamma(y_1,y_2) \right](t-u,u,u), \\
     & \ds = \mu^2 \, \left[ y_1 (\partial_z y_1 +\partial_{\tz} y_1) - \frac{n}{2}\p_z \, y_2\right](t-u,u,u).
\end{array}
\end{equation}
Let us set $$\ty(t,u):=y_1(t-u,u,u)=\frac{1}{\cosh(\mu \, u)}\int_0^u\delta(t-s) \, \cosh(\mu \, s)\, ds.$$We have $$\partial_u \ty(t,u)=[-\partial_t y_1+\partial_z y_1+\partial_{\tz} y_1](t-u,u,u),$$thus,
\begin{equation}\label{eq:y1a}
    [y_1 (\partial_z y_1 +\partial_{\tz} y_1)](t-u,u,u)=[\ty \, (\p_u \ty + \p_t \ty)](t,u).
\end{equation}
Using \eqref{eq:sys_y}, and since $y_2(u)=\tanh(\mu \, u)/\mu$, we get:
\begin{equation}\label{eq:duy1}
  \partial_u \ty(t,u)=\delta(t-u)-\mu \, \ty(t,u) \, \tanh(\mu\, u),
\end{equation}
and differentiating $\ty$ with respect to $t$ and integrating by parts, we get:
\begin{equation}\label{eq:dty1}
    \begin{array}{rl}
     \p_t\ty(t,u)= &\ds \frac{1}{\cosh(\mu \, u)}\int_0^u\delta'(t-s) \, \cosh(\mu \, s)\, ds,\\
    = & \ds \frac{1}{\cosh(\mu \, u)}\left[\delta(t)-\delta(t-u)\,\cosh(\mu\, u)+ \mu \int_0^u\delta(t-s) \, \sinh(\mu \, s)\, ds\right].\\
\end{array}
\end{equation}
Combining \eqref{eq:y1a}, \eqref{eq:duy1} and \eqref{eq:dty1}, we get:
\begin{multline}\label{eq:y1b}
   [y_1 (\partial_z y_1 +\partial_{\tz} y_1)](t-u,u,u)\\=\ty \, \left[\frac{\delta(t)}{\cosh(\mu \, u)}-\mu \, \ty \, \tanh(\mu\, u)+\mu\,\int_0^u\delta(t-s) \, \frac{\sinh(\mu \, s)}{\cosh(\mu \, u)}\, ds\right].
\end{multline}
Next, integrating by parts and using standard trigonometric formulas, we note that:
\begin{equation}\label{eq:y1_int}
    \begin{array}{rl}
     \ds \int_0^t \frac{1}{\cosh(\mu \, u)} \, \ty(t,u) \, du& = \ds \int_0^t \frac{1}{\cosh^2(\mu \, u)} \, \int_0^u\delta(t-s) \, \cosh(\mu \, s)\, ds \, du,\\
     & =  \ds \frac{ 1}{\mu} \left[\tanh(\mu\, t) \int_0^t\delta(t-u)\cosh(\mu\,u )  du -\int_0^t\delta(t-u)\sinh(\mu\,u )  du  \right], \\
     & =  \ds  \frac{1}{\mu^2} \, H_\delta(t),
\end{array}
\end{equation}
with $$H_\delta(t):=\mu \, \int_0^t\delta(u) \, \frac{\sinh(\mu \, u)}{\cosh(\mu \, t)}\,  du.$$Similarly,
\begin{equation}\label{eq:y1_int2}
   \int_0^t \ty(t,u) \, \left[-\mu \, \ty(t,u) \, \tanh(\mu\, u)+\mu\,\int_0^u\delta(t-s) \, \frac{\sinh(\mu \, s)}{\cosh(\mu \, u)}\, ds\right]\, du=-\frac{1}{2\, \mu^2} H_\delta(t)^2.
\end{equation}
Integrating \eqref{eq:y1b} between $0$ and $t$, and using \eqref{eq:y1_int} and \eqref{eq:y1_int2}, we get:
\begin{equation}\label{eq:y1c}
   \int_0^t [y_1 (\partial_z y_1 +\partial_{\tz} y_1)](t-u,u,u) \, du=\frac{1}{\mu^2}\delta(t) \, H_\delta(t)-\frac{1}{2\, \mu^2} H_\delta(t)^2.
\end{equation}
Using \eqref{eq:beta}, we get:
\begin{equation}\label{eq:y1d}
   \int_0^t  (\partial_z \beta+\partial_{\tz} \beta)(t-u,u,u) \, du=\delta(t) \, H_\delta(t)-\frac{1}{2} H_\delta(t)^2-\mu \, \frac{n}{2} \, \tanh(\mu \,t),
\end{equation}
and coming back to \eqref{eq:mb_beta}, this shows that:
\begin{equation}\label{eq:mbar_expl_proof}
   \mb(t)=-\mu \, \frac{n}{2} \, \tanh(\mu \,t) -\frac{1}{2}\lp H_\delta(t)-\delta(t) \rp^2+  R_0'(t),
\end{equation}
with $R_0'(t)= (\p_z Q_0+\p_{\tz} Q_0)(t,t)=(\p_z y_1+\p_{\tz}y_1)(0,t,t) \partial_1 C_0(0,\phi_0(t,t))+\p_z y_2(t) \p_2 C_0(\phi_0(t,t)).$ This proves the formula~\eqref{eq:mbar_cor}.

Lastly,
we note that $$(\p_z y_1+\p_{\tz}y_1)(0,t,t) = -\mu \tanh(\mu \, t)\, \ty(t,t)+ \p_t \ty (t,t),$$and from formula~\eqref{eq:dty1},
\begin{equation}\label{eq:dzy1_0}
    \begin{array}{l}
(\p_z y_1+\p_{\tz}y_1)(0,t,t) = \\
\ds \frac{1}{\cosh(\mu \, t)}\left[-\mu \tanh(\mu \, t) \int_0^t\delta(t-s) \,\cosh(\mu \,s ) \, ds + \mu \, \int_0^t\delta(t-s) \,\sinh(\mu \,s ) \, ds +\delta(t) \right] \\ = \ds \frac{\mu}{\cosh^2(\mu \, t)}\left[ \int_0^t\delta(t-s)\, \lp\sinh(\mu \,s ) \, \cosh(\mu \, t) -\sinh(\mu \,t ) \,\cosh(\mu \,s ) \rp\, ds \right] +\frac{\delta(t)}{\cosh(\mu \, t)}\\
=\ds \frac{\mu}{\cosh^2(\mu \, t)}\left[ \int_0^t\delta(t-s)\, \sinh(\mu \,(s-t)) \, ds \right] +\frac{\delta(t)}{\cosh(\mu \, t)} \\
= \ds \frac{1}{\cosh(\mu \, t)}\lp \delta(t)-H_\delta(t) \rp.
\end{array}
\end{equation}
Finally, this yields:
\begin{equation}\label{eq:Q0prime}
    R_0'(t)=  \frac{1}{\cosh(\mu \, t)}\lp \delta(t)-H_\delta(t) \rp  \partial_1 C_0(0,\phi_0(t,t)) +(1- \tanh^2(\mu \,t ))\,\p_2 C_0(\phi_0(t,t)),
\end{equation}
with $\phi_0(t,t)=(y_1(0,t,t),y_2(t)).$ This concludes the proof of Corollary~\ref{cor:formule_mbar}. $\Box$

\subsection{Proof of Corollary~\ref{cor:formule_variance}}
In order to simplify the computations, we introduce in this section the function
\begin{equation}\label{def:y1hat}
 \hy(t,z):= y_1(t,z,z).
\end{equation}
Notice that:
\begin{equation}\label{eq hy}
 \left\{
 \begin{array}{l}
\hy (t,z) = \ds\int_0^{z} \delta(z+t-s) \, \frac{\cosh(\mu \, s)}{\cosh(\mu \, z)} ds\, \\
\partial_t \hy -\partial_z \hy =-\delta(t)+\mu^2 \,\hy \,y_2,  \\
\hy(t,0)=0, \ \partial_t \hy(t,0)= 0,
 \end{array}
 \right.
\end{equation}
and, finally with $R$ defined by $R(t,z):=Q(t,z,z)$,
\begin{align*}
    R(t,z) = C(t,\hy,y_2).
\end{align*}
Then, from a straightforward computation, we get:
\begin{multline}
    \partial_{zz} R(t,0)=(\partial_z \hy(t,0))^2  \partial_{11} C(t,0,0)+ (\partial_z y_2(0) )^2  \partial_{22} C(t,0,0)  \\
    + 2 \partial_z \hy(t,0)\partial_z y_2(0) \partial_{12} C(t,0,0) +\partial_1 C(t,0,0)  \p_{zz} \hy(t,0) + \partial_2 C(t,0,0) \p_{zz}y_2(0).
\end{multline}
Thanks to \eqref{eq hy},
$$
    \partial_{zz} R(t,0)= \delta(t)^2  \partial_{11} C(t,0,0)+  \partial_{22} C(t,0,0)
    + 2 \delta(t) \partial_{12} C(t,0,0) +\partial_1 C(t,0,0)  \delta'(t) .$$Now coming back to the formula for the variance we established thanks to the CGF in \eqref{eq:var_cgf},
\begin{align*}
  V_m(t) =  \p_{zz} R(t,0) -\partial_1 C(t,0,0)  \delta'(t).
\end{align*}
Next, from \Li Theorem~\ref{th:formule_C}\Bk, differentiating with respect to $\tz$, we get:
\begin{align*}
    \p_{\tz} Q(t,0,0) = \p_1 C(t,0,0) \p_{\tz} y_1(t,0,0).
\end{align*}
Since $ \p_{\tz} y_1(t,0,0) = -\cosh(\mu t), $
we obtain the expression of the variance given in \Li Corollary~\ref{cor:formule_variance}.\Bk

\section*{Appendix A: relationship between the equations \eqref{eq:main_PDE_moving_env} and \eqref{eq:main_form_n}}

Consider $q(t,\x)$ the solution of \eqref{eq:main_PDE_moving_env} and set  $r(t,\x)=r_{max}+m(t,\x)$, with $r_{max}>0$ and $m(t,\x)$ defined by \eqref{eq:m(t,x)}. Define the `total population' at time $t$ as the solution of
\begin{equation}\label{eq:rhot}
    \rho'(t)= \rho(t) \, (\rb(t) -\rho(t)),
\end{equation}
with$$\rb(t):=\int_{\R^n} r (t,\x) \, q(t,\x) \, d\x,$$the mean growth rate in the population at time $t$. Then, setting $n(t,\x):= \rho(t) \, q(t,\x)$ (the population density) we observe that
$$
\partial_t n (t, \x) = \frac{\mu^2}{2} \Delta n+  n(t, \x) \, (r (t,\x)-\rho(t)), \  t >0, \  \x \in \R^n.
$$
Thus, if $\mb(t)$ has a limit $\mb(\infty)$ as $t\to +\infty$ (e.g., in the case of a linearly or sublinearly moving optimum, see Propositions~\ref{prop:linear} and~\ref{prop:nonlinear} (i)) the population size $\rho(t)$ converges to $r_{max}+\mb(\infty)$. In particular, large-time persistence is equivalent to $r_{max}+\mb(\infty)>0$.

\Li
\noindent Conversely, consider $n(t,\x)$ a positive solution of~\eqref{eq:main_form_n} with $\rho(t)=\int_{\Omega}n(t,\x)\, d \x>0$ and $r(t,\x)\le r_{max}$, with $r_{max}$ a positive constant. Define $q(t,\x)=n(t,\x)/\rho(t)$. Integrating the equation \eqref{eq:main_form_n} over $\Omega\subseteq\R^n$, and provided that $\int_{\Omega}\Delta n=0$ (which means that the mutations do not change the total mass) we note that $\rho(t)$ satisfies~\eqref{eq:rhot}. Thus,
$$
\partial_t q (t, \x) =\partial_t n (t, \x)/\rho(t)-q(t,\x)\,\rho'(t)/\rho(t)=  \frac{\mu^2}{2} \Delta q+  q(t, \x) \, (r (t,\x)-\rb(t)), \  t >0, \  \x \in \Omega\subseteq\R^n,
$$
and finally, $q$ solves an equation of the form \eqref{eq:main_PDE_moving_env}, with $m(t,\x)=r(t,\x)-r_{max}.$
\Bk

\section*{Appendix B: Skewness}

As the fitness satisfies \eqref{eq:m(t,x)_b}, $m(t,\x)=\delta(t) \, \mm_1(\x)+\mm_2(\x)-\delta(t)^2/2$, the third central moment $\M_3(t)$ of the variable $m(t,\x)$ is equal to the third central moment of $\delta(t) \, \mm_1(\x)+\mm_2(\x)$. Let us define$$S(t,z):=C(t,\delta(t) \, z, z), \ t\ge 0, \ z\ge 0.$$The function $S(t,z)$ satisfies:
$$S(t,z)=\ln\lp \int_{\R \times \R_-} p(t,m_1,m_2) e^{z(\delta(t)\, m_1 + m_2)} \, dm_1 \, dm_2\rp,$$and therefore corresponds to the cumulant generating function of $\delta(t) \, \mm_1(\x)+\mm_2(\x)$. Its third central moment is therefore given by:
$\M_3(t)=\partial_{zzz}S(t,0).$
The skewness of the distribution of fitness is therefore given by the formula:
\begin{align*}
   \hbox{Skew}_m(t)= \frac{\p_{zzz} S(t,0)}{V_m(t)^\frac 32}.
\end{align*}
By differentiating, one finds that$$
    \p_{zzz} S(t,0) = \delta(t)^3 \p_{111} C(t,0,0) + 3 \delta(t)^2 \p_{112} C(t,0,0)+ 3 \delta(t) \p_{122} C(t,0,0) + \p_{222} C(t,0,0).$$Moreover, since $\p_{zz} y_2 (0)=0 $,
\begin{align*}
    \partial_{zzz} R(t,0)=&\delta(t)^3 \partial_{111} C(t,0,0)+ 3 \delta(t) \delta'(t) \partial_{11} C(t,0,0)+ 3 \delta(t)^2 \partial_{112} C(t,0,0) \\ & +  \partial_{222} C(t,0,0) + 3\delta(t) \partial_{221} C(t,0,0) + 3 \delta'(t) \partial_{12} C(t,0,0) \\
    &+ \partial_{zzz} \hy(t,0) \p_1 C(t,0,0)
    +\partial_{zzz} y_2(0) \p_2 C(t,0,0),
\end{align*}
with $\hy$ defined by \eqref{def:y1hat}. Therefore, the skewness is equal to:
\begin{multline*}
    \hbox{Skew}_m(t)=\frac{\p_{zzz} R(t,0)}{V_m(t)^{3/2}} \\
    -\frac{ 3 \delta(t) \delta'(t) \partial_{11} C(t,0,0)+ 3 \delta'(t) \partial_{12} C(t,0,0) + \partial_{zzz} \hy(t,0) \p_1 C(t,0,0)
     + \partial_{zzz} y_2(0) \p_2 C(t,0,0)}{V_m(t)^{3/2}}.
\end{multline*}
Next, straightforward computations yield
$$
\begin{array}{rcl}
   \p_{\tz \tz} Q(t, 0, 0)  &=&   \p_{11} C(t,0,0)\cosh(\mu t)^2, \\
    \p_{z \tz}Q(t, 0, 0)& = &- ( \delta(t)+\cosh(\mu\,t)) \p_{11} C(t,0,0)\cosh(\mu t)- \p_{12}C(t,0,0)\cosh(\mu t)\\
    & & -\mu \sinh(\mu t) \p_1 C(t,0,0) ,\\
    \p_{zzz} y_2(0)& =& -2 \mu^2, \\
    \p_{zzz}\hy(t,0) &=& -2 \mu ^2 \delta(t) + \delta''(t).
\end{array}
$$Finally, one finds that
\begin{multline}\label{form skew}
    \hbox{Skew}_m(t)=\frac{1}{V_m(t)^{3/2}} \left[ \p_{zzz} R(t,0) + 2 \mu^2 \p_z R(t,0)
     + \delta''(t)  \frac{ \partial_{\tz } Q(t,0,0) } {\cosh(\mu t)}  \right. \\ \left. +  3 \delta'(t) \frac{  \partial_{z \tz} Q(t,0,0) +  \p_{\tz \tz} Q(t,0,0)-\mu \tanh(\mu t) \p_{\tz} Q(t,0,0)}{\cosh(\mu t)}
      \right].
\end{multline}

\section*{Acknowledgements}
This work was supported by the French Agence Nationale de la Recherche (ANR-18-CE45-0019 ``RESISTE"). The authors thank the reviewers for valuable comments and suggestions.

\bibliographystyle{plain}
%\bibliography{biblio_lionel_jab_drive}

\end{document}